\documentclass[12pt]{amsart}
\usepackage{amsmath}
\usepackage{amssymb}
\usepackage{latexsym}

\newcommand{\Zint}{\mathbb {Z}}    
     
\newcommand{\Cplx}{\mathbb {C}}     

\newcommand{\halmos}{\rule{5pt}{5pt}}

\numberwithin{equation}{section}

\newtheorem{prop}{\bf Proposition}[section]

\newtheorem{thm}[prop]{\bf Theorem}

\newtheorem{cor}[prop]{\bf Corollary}

\begin{document}

\title[Heun equation {\rm V}]
{The Heun equation and the Calogero-Moser-Sutherland system {\rm V}: generalized Darboux transformations}
\author{Kouichi Takemura}
\address{Department of Mathematical Sciences, Yokohama City University, 22-2 Seto, Kanazawa-ku, Yokohama 236-0027, Japan.}
\email{takemura@yokohama-cu.ac.jp}


\begin{abstract}
We obtain isomonodromic transformations for Heun's equation by generalizing Darboux transformation, and we find pairs and triplets of Heun's equation which have the same monodromy structure. By composing generalized Darboux transformations, we establish a new construction of the commuting operator which ensures finite-gap property. As an application, we prove conjectures in part III.
\end{abstract}

\maketitle

\section{Introduction}

It was shown in \cite{Tak4} that some pairs of Schr\"odinger operators are isomonodromic. Set
\begin{align}
& H_1= -\frac{d^2}{dx^2} + 6\wp (x), \label{Ino1} \\
& H_2= -\frac{d^2}{dx^2} + 2\wp (x) +2\wp (x+\omega_1) +2\wp (x+\omega_2) , \label{Ino2}
\end{align}
where $\wp (x)$ is the Weierstrass $\wp$-function with periods $(2\omega _1,2\omega _3 )$ and $\omega _2=-\omega _1-\omega _3$.
Now we consider eigenfunctions of $H_1$ (resp. $H_2$) with the eigenvalue $E$.
Set
\begin{align*}
 \Xi _1 (x,E)= & 9\wp(x)^2+3E\wp(x)+E^2-9g_2/4, \\
 \Xi _2 (x,E)= & (E-3e_3)\wp(x)+(E-3e_2)\wp(x +\omega _1) \\
 & +(E-3e_1)\wp(x +\omega _2) +E^2-3g_2/2, \nonumber \\
 Q(E)= & (E^2-3g_2)(E-3e_1)(E-3e_2)(E-3e_3) , \\
 \Lambda _k (x,E)  = & \sqrt{\Xi _k(x,E)}\exp \int \frac{\sqrt{-Q(E)}dx}{\Xi _k(x,E)} , \; \; \; (k=1,2), 
\end{align*}
where $e_i =\wp (\omega _i)$ $(i=1,2,3)$ and $g_2=-4(e_1e_2+e_2e_3+e_3e_1)$. Then it was shown that the functions $\Lambda _1(x,E)$ and $\Lambda _1(-x,E)$ (resp. $\Lambda _2(x,E)$ and $\Lambda _2(-x,E)$) are eigenfunctions of $H_1$ (resp. $H_2$) with the eigenvalue $E$, and they satisfy
\begin{equation*} 
\Lambda _k (\pm (x+2\omega _i),E)= \Lambda _k (x,E) \exp \left( \mp \frac{1}{2} \int_{\sqrt{3g_2}}^{E}\frac{ -6\tilde{E} \eta _i +(2\tilde{E}^2-3g_2)\omega _i }{\sqrt{-Q(\tilde{E})}} d\tilde{E}\right) ,
\end{equation*}
for $k=1,2$ and $i=1,2,3$, where $\eta _i =\zeta (\omega _i)$ and $\zeta (x)$ is the Weierstrass zeta function. Hence the monodromy of eigenfunctions of $H_1$ with the eigenvalue $E$ coincides with that of eigenfunctions of $H_2$.

In this paper, we investigate this phenomena by Darboux transformation and generalized Darboux transformation.
Let $\phi _0 (x)$ be an eigenfunction of the operator $H= -d^2/dx^2 + q(x)$ with an eigenvalue $E_0$, i.e. 
\begin{equation*}
\left( -\frac{d^2}{dx^2} + q(x) \right) \phi _0(x) = E_0 \phi _0(x).
\end{equation*}
For this case, the potential $q(x)$ is written as $q(x)=(\phi _0'(x)/\phi _0(x))'+(\phi _0'(x)/\phi _0(x))^2+E_0$. Set $L= d/dx - \phi _0'(x)/\phi _0(x)$ and $\tilde {H}=  -d^2/dx^2 +q(x) - 2(\phi _0'(x)/\phi _0(x))'$.
Then we have
\begin{equation*}
\tilde{H}L=LH.
\end{equation*}
Hence, if $\phi (x)$ is an eigenfunction of the operator $H$ with the eigenvalue $E$, then $L \phi (x)$ is an eigenfunction of the operator $\tilde{H}$ with the eigenvalue $E$. This transformation is called the Darboux transformation.
We generalize the operator $L$ to be the differential operator of higher order, and we call it the generalized  Darboux transformation.

The Schr\"odinger operator we consider in this paper is the Hamiltonian of the $BC_1$ Inozemtsev model, which is written as 
\begin{equation}
H^{(l_0,l_1,l_2,l_3)}= -\frac{d^2}{dx^2} + \sum_{i=0}^3 l_i(l_i+1)\wp (x+\omega_i),
\label{Ino}
\end{equation}
where $\omega _0 =0$. The potential of this operator is called the Treibich-Verdier potential, because Treibich and Verdier \cite{TV} found and showed that, if $l_i \in \Zint _{\geq 0}$ for all $i \in \{ 0,1,2,3 \}$, then it is an algebro-geometric finite-gap potential. For further results on this subject, see \cite{GW,Smi,Tak1,Tak3,Tre}. The algebro-geometric finite-gap property cause the possibility for calculation of eigenfunction and monodromy of the operator $H^{(l_0,l_1,l_2,l_3)}$.

Let $f(x)$ be an eigenfunction of the operator $H^{(l_0,l_1,l_2,l_3)}$ with the eigenvalue $E$, namely,
\begin{equation}
\left( -\frac{d^2}{dx^2} + \sum_{i=0}^3 l_i(l_i+1)\wp (x+\omega_i) \right) f(x)= Ef(x).
\label{Heq}
\end{equation}
Then this equation is an elliptic representation of Heun's equation. Here Heun's equation is the standard canonical form of a Fuchsian equation with four singularities (see \cite{Ron}). Thus, solving Heun's equation is equivalent to studying eigenvalues and eigenfunctions of the Hamiltonian of the $BC_1$ Inozemtsev model.

We now describe the main result of this paper.
Let $\alpha_i$ be a number such that $\alpha_i= -l_i$ or $\alpha_i= l_i+1$ for each $i\in \{ 0,1,2,3\} $.
Set $d=-\sum_{i=0}^3 \alpha_i /2$ and assume $d\in \Zint_{\geq 0}$.
Then there exists a differential operator $L$ of order $d+1$ which satisfies
\begin{equation*}
H^{(\alpha _0 +d,\alpha _1 +d,\alpha _2 +d,\alpha _3 +d)}  L=L H^{(l_0,l_1,l_2,l_3)}.
\end{equation*}
Note that Khare and Sukhatme \cite{KS} essentially established this result for the case $d=0$, that is the case of original Darboux transformation. 
It follows immediately that, if $\phi (x)$ is an eigenfunction of the operator $H^{(l_0,l_1,l_2,l_3)}$ with an eigenvalue $E$, then  $L\phi (x)$ is an eigenfunction of the operator $H^{(\alpha _0 +d,\alpha _1 +d,\alpha _2 +d,\alpha _3 +d)}$ with the eigenvalue $E$.
Since all coefficients of the operator $L$ with respect to the differential $(d/dx)^k$ $(k=0,\dots ,d+1)$ is shown to be doubly-periodic, the operator $L$ preserves the data of monodromy.
Hence the operators $H^{(l_0,l_1,l_2,l_3)}$ and $H^{(\alpha _0 +d,\alpha _1 +d,\alpha _2 +d,\alpha _3 +d)}$ are isomonodromic, and isospectral, because boundary condition for spectral problem is characterized by monodromy. Note that the condition $d\in \Zint_{\geq 0}$ corresponds to quasi-solvability of the operator $H^{(l_0,l_1,l_2,l_3)}$.

For the case that $l_0, l_1, l_2, l_3$ are all integers, there exists an operator $H^{(\tilde{l}_0,\tilde{l}_1,\tilde{l}_2,\tilde{l}_3)}$ such that the pair $H^{(l_0,l_1,l_2,l_3)}$ and $H^{(\tilde{l}_0,\tilde{l}_1,\tilde{l}_2,\tilde{l}_3)}$ is connected by isomonodromic transformation. In some cases, they are self-dual. For example, the operator $H_1(=H^{(2,0,0,0)})$ in Eq.(\ref{Ino1}) is connected to the operator $H_2(=H^{(1,1,1,0)}$) in Eq.(\ref{Ino2}) by the transformation $L=d/dx -\wp ' (x) /(2(\wp(x)-e_1)) -\wp ' (x)/(2(\wp (x)-e_2 ))$, i.e. we have
\begin{equation*}
H^{(1,1,1,0)} L=L H^{(2,0,0,0)}.
\end{equation*}
For the case that $l_0 +1/2, l_1 +1/2, l_2 +1/2, l_3  +1/2$ are all integers, there exists two operators $H^{(l_0 ^{(1)},l_1^{(1)},l_2^{(1)},l_3^{(1)})}$ and $H^{(l_0 ^{(2)},l_1^{(2)},l_2^{(2)},l_3^{(2)})}$ such that the triplet $H^{(l_0,l_1,l_2,l_3)}$, $H^{(l_0 ^{(1)},l_1^{(1)},l_2^{(1)},l_3^{(1)})}$ and $H^{(l_0 ^{(2)},l_1^{(2)},l_2^{(2)},l_3^{(2)})}$ is connected by isomonodromic transformations.

In the paper \cite{Tak3}, finite-gap property of the operator $H^{(l_0,l_1,l_2,l_3)}$ for the case $l_0, l_1, l_2, l_3 \in \Zint _{\geq 0}$ is studied (see also \cite{TV,Tre,GW,Smi}). Especially, a differential operator $A$ of odd order which commutes with $H^{(l_0,l_1,l_2,l_3)}$ is constructed.
In this paper, we propose a new method for construction of the commuting operator by composing four generalized Darboux transformations. Note that each generalized Darboux transformation is written explicitly.
To show that the commuting operator constructed by composing four generalized Darboux transformations coincides with the one defined in \cite{Tak3}, we need a discussion which will be done in section \ref{sec:fin}.
As an application, we prove conjectures in \cite{Tak3}.
Namely, we establish that the polynomial defined by quasi-solvability coincides with the one defined by using the doubly-periodic function $\Xi (x,E)$ which is written as a product of two eigenfunctions. We also prove that the commuting operator is characterized by annihilating the spaces of quasi-solvability.
As for the isomonodromic pair $H^{(l_0,l_1,l_2,l_3)}$ and $H^{(\tilde{l}_0,\tilde{l}_1,\tilde{l}_2,\tilde{l}_3)}$, it is shown that some functions related to the monodromy of $H^{(l_0,l_1,l_2,l_3)}$ coincide with that of $H^{(\tilde{l}_0,\tilde{l}_1,\tilde{l}_2,\tilde{l}_3)}$.

This paper is organized as follows. In section \ref{sec:QSD}, we review the connection between quasi-solvability and generalized Darboux transformation which was essentially done in \cite{AST}.
In section \ref{sec:Heun}, we consider generalized Darboux transformations for the case of Heun's equation.
In section \ref{sec:int} (resp. section \ref{sec:hint}), we investigate isomonodromic transformations for the case $l_0, l_1, l_2, l_3 \in \Zint$ (resp. $l_0, l_1, l_2, l_3 \in \Zint +1/2$).
In section \ref{sec:fin}, we construct a differential operator of odd order which commutes with $H^{(l_0,l_1,l_2,l_3)}$ and investigate it.

\section{Quasi-solvability and generalized Darboux transformation} \label{sec:QSD}

We review the relationship between the quasi-solvability and the generalized Darboux transformation.

We set $H=-d^2/dx^2 +q(x)$. Let $n$ be a positive integer. If the operator $H$ preserve a $n$-dimensional space $U$ of functions, then the operator $H$ is called quasi-solvable. Then there exists a basis $\langle f_1 (x), \dots , f_{n} (x) \rangle $ of the invariant space such that $H f_j (x)= \sum _i a_{i,j} f _i(x)$ for some constants $a_{i,j}$ $(1\leq i,j \leq n)$. Let $P_{H,U} (t)$ be the characteristic polynomial of the operator $H$ on the space $U$. Then the set $\{ E| P_{H,U} (E)=0\}$ coincides with the set of eigenvalues of of the operator $H$ on the space $U$. Then the model is partially solved, and this is an origin of ``quasi-solvability''.
For the space $U$, there exists a monic differential operator of order $n$
\begin{equation}
L= \left( \frac{d}{dx} \right) ^n +\sum _{i=1}^n c_i(x) \left( \frac{d}{dx} \right) ^{n-i},
\label{eq:L}
\end{equation}
such that 
 $L f(x) =0$ for all $f(x) \in U$. It is determined uniquely and written as 
\begin{equation*}
L=
\begin{scriptsize}
\left| \left(
\begin{array}{cccc}
f_1(x) & \frac{d}{dx}f_1(x) & \dots &  \left(\frac{d}{dx} \right)^{n-1} f_1(x)\\
f_2(x) & \frac{d}{dx}f_2(x) & \dots &  \left(\frac{d}{dx} \right)^{n-1} f_2(x)\\
\vdots & & \vdots & \vdots \\
f_n(x) & \frac{d}{dx}f_n(x) & \dots &  \left(\frac{d}{dx} \right)^{n-1} f_n(x)
\end{array} \right) \right| ^{-1}
 \left| \left(
\begin{array}{cccc}
f_1(x) & \frac{d}{dx}f_1(x) & \dots &  \left(\frac{d}{dx} \right)^{n} f_1(x)\\
f_2(x) & \frac{d}{dx}f_2(x) & \dots &  \left(\frac{d}{dx} \right)^{n} f_2(x)\\
\vdots & & \vdots & \vdots \\
f_n(x) & \frac{d}{dx}f_n(x) & \dots &  \left(\frac{d}{dx} \right)^{n} f_n(x)\\
1 & \frac{d}{dx} & \dots &  \left(\frac{d}{dx} \right)^{n} 
\end{array} \right) \right| .
\end{scriptsize}
\end{equation*}

\begin{prop} \label{prop:GDT} (c.f. \cite{AST})
Assume that the operator $H=-d^2/dx^2 +q(x)$ preserve a $n$-dimensional space $U$ of functions. Let $L$ be the differential operator written as Eq.(\ref{eq:L}) which annihilate the functions in $U$. Set $\tilde{H} = -d^2/dx^2 +q(x)+2c'_1(x)$. Then we have
\begin{equation*}
\tilde{H}L=LH.
\end{equation*}
\end{prop}
\begin{proof}
By a direct calculation, it follows that the order of the differential operator $\tilde{H}L-LH$ is at most $n-1$. Assume that $\tilde{H}L-LH\neq 0$. We denote the order by $k$. Then the dimension of solutions to the differential equation $(\tilde{H}L-LH)f(x)=0$ is $k$. Let $g(x) \in U$. Since the operator $H$ preserve the space $U$, we have $H g(x) \in U$. By definition of the space $U$, we have $L g(x)=0$ and $LHg(x)=0$. Hence we have $(\tilde{H}L-LH)g(x)=0$ for $g(x) \in U$, but it contradicts to that the dimension of solutions is $k(\leq n-1)$. Therefore we obtain $\tilde{H}L=LH$.
\end{proof}
We consider the case $n=1$, Let $\phi _0 (x)$ be a non-zero function in $U$. Then $U= \Cplx \phi _0 (x)$, and the operator which annihilate $\phi _0 (x)$ is written as $L=d/dx - \phi' _0 (x)/\phi _0 (x)$. The operator $\tilde {H}$ is written as $\tilde {H}=  H - 2(\phi _0'(x)/\phi _0(x))'$. Hence the proposition reproduce the Darboux transformation. In this sense, the transformation in the proposition may be called the generalized Darboux transformation.

\section{Generalized Darboux transformation for Heun's equation} \label{sec:Heun}

In this section, we apply Proposition \ref{prop:GDT} for Heun's equation. For this purpose, we recall quasi-solvability of Heun's equation.
\begin{prop} \label{findim} \cite[Proposition 5.1]{Tak2}
Let $\alpha _i$ be a number such that $\alpha _i= -l_i$ or $\alpha _i= l_i+1$ for each $i\in \{ 0,1,2,3\} $.
Set $d=-\sum_{i=0}^3 \alpha _i /2$ and assume $d\in \Zint_{\geq 0}$.
Let $V_{\alpha _0, \alpha _1, \alpha _2, \alpha _3}$ be the $d+1$-dimensional space spanned by 
\begin{equation}
\left\{ (\wp (x) -e_1)^{\alpha _1/2}(\wp (x) -e_2)^{\alpha _2/2}(\wp (x) -e_3)^{\alpha _3/2} \wp(x)^n\right\} _{n=0, \dots ,d}. \label{eq:Vqsol}
\end{equation}
Then the operator $H^{(l_0,l_1,l_2, l_3)}$ (see Eq.(\ref{Ino})) preserves the space $V_{\alpha _0, \alpha _1, \alpha _2, \alpha _3}$.
\end{prop}
Set $z=\wp(x)$, $\widehat{\Phi}(z)=(z-e_1)^{\alpha _1/2}(z-e_2)^{\alpha _2/2}(z-e_3)^{\alpha _3/2}$, and $\widehat{H}^{(l_0,l_1,l_2, l_3)}= \widehat{\Phi}(z)^{-1} \circ H^{(l_0,l_1,l_2, l_3)} \circ\widehat{\Phi}(z)$. Proposition \ref{findim} is proved by showing that the operator $\widehat{H}^{(l_0,l_1,l_2, l_3)}$ preserve the space spanned by $(z-e_2)^r$ $(r=0,\dots ,d)$. For details, see the proof of \cite[Proposition 5.1]{Tak2}.

Now we calculate the differential operator which annihilate the space $V_{\alpha _0, \alpha _1, \alpha _2, \alpha _3}$.
\begin{prop} \label{prop:L0123}
The monic differential operator of order $d+1$ which annihilate the space $V_{\alpha _0, \alpha _1, \alpha _2, \alpha _3}$ is written as
\begin{equation}
L_{\alpha _0, \alpha _1, \alpha _2, \alpha _3}= \wp '(x)^{d+1} \widehat{\Phi}(\wp (x)) \circ \left(\frac{1}{\wp' (x)}\frac{d}{dx} \right)^{d+1} \circ \widehat{\Phi}(\wp (x)) ^{-1}.
\label{op:L}
\end{equation}
We write the operator $L_{\alpha _0, \alpha _1, \alpha _2, \alpha _3}$ as
\begin{equation}
L_{\alpha _0, \alpha _1, \alpha _2, \alpha _3}= \left( \frac{d}{dx} \right) ^{d+1} +\sum _{i=1}^{d+1} c_i(x) \left( \frac{d}{dx} \right) ^{d+1-i}.
\label{eq:L0123}
\end{equation}
Then  
\begin{equation}
c_1 (x)  = -\frac{d+1}{4} \left( \sum _{i=1}^3 \frac{2\alpha _i+d}{\wp (x)-e_i} \right)\wp '(x). \label{eq:a1}
\end{equation}
If $i$ is even (resp. odd), then $c_{i} (x)$ is expressed as $c_{i} (x)=R_{i} (\wp (x))$ (resp. $c_{i} (x)=R_{i} (\wp (x))\wp '(x)$), where $R_{i} (z)$ is a rational function in $z$.
\end{prop}
\begin{proof}
It is trivial that the operator $(d/dz)^{d+1}$ annihilate the space spanned by $z^r$ $(r=0,\dots ,d)$, and the operator $ \widehat{\Phi}(z) \circ (d/dz)^{d+1} \circ  \widehat{\Phi}(z) ^{-1}$ annihilate the space spanned by $\widehat{\Phi}(z)  z^r$ $(r=0,\dots ,d)$. Write 
\begin{equation*}
\widehat{\Phi}(z)  \circ \left( \frac{d}{dz} \right)^{d+1} \circ  \widehat{\Phi}(z) ^{-1}= \left( \frac{d}{dz} \right) ^{d+1} +\sum _{i=1}^{d+1} \hat{c}_i(z) \left( \frac{d}{dz} \right) ^{d+1-i}.
\end{equation*}
Then $\hat{c}_1(z)= -\sum _{i=1}^3 ((d+1)\alpha _i)/(2(z-e_i))$, and  all coefficients $\hat{c}_i(z)$ are rational functions in $z$.
By the transformation $z=\wp (x)$, the monic operator 
\begin{align}
L_{\alpha _0, \alpha _1, \alpha _2, \alpha _3}& =  \wp '(x)^{d+1} \widehat{\Phi}(\wp (x)) \circ \left(\frac{1}{\wp' (x)}\frac{d}{dx} \right)^{d+1} \circ \widehat{\Phi}(\wp (x)) ^{-1} \label{eq:Latrans} \\
& = \wp '(x)^{d+1} \left( \frac{1}{\wp' (x)}\frac{d}{dx} \right) ^{d+1} +\sum _{i=1}^{d+1} \hat{c}_i(\wp (x) ) \wp '(x)^{d+1}  \left( \frac{1}{\wp' (x)}\frac{d}{dx} \right) ^{d+1-i} \nonumber
\end{align}
annihilate the space $V_{\alpha _0, \alpha _1, \alpha _2, \alpha _3}$. Write $\left( (1/\wp' (x))(d/dx) \right) ^{i}= \sum_{j=0}^i \hat{b} _j(x) (d/dx) ^{j}$. If $j$ is even (resp. odd), then $\hat{b} _j(x)$ is expressed as $r _j(\wp (x)) $ (resp. $r _j(\wp (x)) \wp '(x)$), where $r _j(z)$ is a rational function in $z$.
Combining with the relation $\wp ' (x)^2 =4(\wp (x)-e_1)(\wp (x)-e_2)(\wp (x)-e_3)$, it follows that, if $i$ is even (resp. odd), then $c_{i} (x)$ is expressed as $c_{i} (x)=R_{i} (\wp (x))$ (resp. $c_{i} (x)=R_{i} (\wp (x))\wp '(x)$), where $R_{i} (z)$ is a rational function in $z$.
We now calculate $c_1(x)$. By Eqs.(\ref{eq:Latrans}, \ref{eq:Leg}), we have
\begin{align*}
& c_1 (x)= -\frac{(d+1)d}{2} \frac{\wp '' (x)}{\wp ' (x)}+ \hat{c}_1(\wp (x) ) \wp '(x)  = -\frac{d+1}{4} \left( \sum _{i=1}^3 \frac{2\alpha _i+d}{\wp (x)-e_i} \right)\wp '(x) . 
\end{align*}
\end{proof}

\begin{thm} \label{thm:HL0123L0123H}
Let $\alpha_i$ be a number such that $\alpha_i= -l_i$ or $\alpha_i= l_i+1$ for each $i\in \{ 0,1,2,3\} $. Set $d=-\sum_{i=0}^3 \alpha_i /2$ and assume $d\in \Zint_{\geq 0}$. Let $L_{\alpha _0, \alpha _1, \alpha _2, \alpha _3}$ be the operator defined in Proposition \ref{prop:L0123}. Then we have
\begin{equation}
H^{(\alpha _0 +d,\alpha _1 +d,\alpha _2 +d,\alpha _3 +d)}  L_{\alpha _0, \alpha _1, \alpha _2, \alpha _3}=L _{\alpha _0, \alpha _1, \alpha _2, \alpha _3} H^{(l_0,l_1,l_2,l_3)}. \label{eq:HL0123L0123H}
\end{equation}
\end{thm}
\begin{proof}
It follows from Propositions \ref{prop:GDT} and \ref{prop:L0123} that
\begin{equation*}
(H^{(l_0,l_1,l_2,l_3)} + 2c' _1(x))  L_{\alpha _0, \alpha _1, \alpha _2, \alpha _3}=L _{\alpha _0, \alpha _1, \alpha _2, \alpha _3} H^{(l_0,l_1,l_2,l_3)},
\end{equation*}
where $c_1 (x)$ is defined in Eq.(\ref{eq:a1}).
By Eq.(\ref{eq:Leg}) we have
\begin{align*}
2 c ' _1 (x) & = \frac{d+1}{2} \left( \sum _{i=1}^3 \frac{2\alpha _i+d}{(\wp (x)-e_i)^2} \right)\wp '(x)^2 -\frac{d+1}{2} \left( \sum _{i=1}^3 \frac{2\alpha _i+d}{\wp (x)-e_i} \right)\wp ''(x) \\
& = -(d+1)(2(\alpha _1 +\alpha _2 +\alpha _3) +3d) \wp (x) + \sum _{i=1}^3 (d+1)(2\alpha _1 +d) \wp (x+\omega_i). \nonumber 
\end{align*}
Since $\alpha_i= -l_i$ or $\alpha_i= l_i+1$, we have $l_i (l_i +1)= - \alpha _i (- \alpha _i +1)$. Hence we obtain $H^{(l_0,l_1,l_2,l_3)} + 2c' _1(x) = H^{(\alpha _0 -1, \alpha _1 -1,\alpha _2 -1,\alpha _3 -1)} + 2c' _1(x)= H^{(\alpha _0 +d,\alpha _1 +d,\alpha _2 +d,\alpha _3 +d)}$ and Eq.(\ref{eq:HL0123L0123H}).
\end{proof}
We consider the converse relation to Eq.(\ref{eq:HL0123L0123H}).
The operator $H^{(\alpha _0 +d,\alpha _1 +d,\alpha _2 +d,\alpha _3 +d)}$ preserve the $d+1$-dimensional space $V_{-\alpha _0-d, -\alpha _1-d, -\alpha _2-d, -\alpha _3-d}$, because $-\alpha _i-d \in \{ -(\alpha _i +d), \alpha _i +d +1 \}$ and $-\sum _{i=0}^3 (-\alpha _i-d ) /2 = d$.
\begin{prop} \label{prop:HL0dL0dH}
(i) We have
\begin{equation}
H^{(l_0,l_1,l_2,l_3)}  L_{-\alpha _0-d, -\alpha _1-d, -\alpha _2-d, -\alpha _3-d}=L_{-\alpha _0-d, -\alpha _1-d, -\alpha _2-d, -\alpha _3-d} H^{(\alpha _0 +d,\alpha _1 +d,\alpha _2 +d,\alpha _3 +d)} . \label{eq:HL0dL0dH}
\end{equation}
(ii) The characteristic polynomial of the operator $H^{(l_0,l_1,l_2,l_3)}$ on the space $V_{\alpha _0, \alpha _1, \alpha _2, \alpha _3}$ coincides with that of the operator $H^{(\alpha _0 +d,\alpha _1 +d,\alpha _2 +d,\alpha _3 +d)}$ on the space $V_{-\alpha _0-d, -\alpha _1-d, -\alpha _2-d, -\alpha _3-d}$.
\end{prop}
\begin{proof}
Set $\tilde{l} _i = \alpha _i +d$ $(i=0,1,2,3)$. Then $-\sum _{i=0}^3 \tilde{l} _i /2 = d$. By Theorem \ref{thm:HL0123L0123H}, we have
\begin{equation*}
H^{(-\tilde{l}_0+d,-\tilde{l}_1+d,-\tilde{l}_2+d,-\tilde{l}_3+d)}  L_{-\tilde{l}_0, -\tilde{l}_1, -\tilde{l}_2, -\tilde{l}_3}=L_{-\tilde{l}_0, -\tilde{l}_1, -\tilde{l}_2, -\tilde{l}_3} H^{(\tilde{l}_0, \tilde{l}_1, \tilde{l}_2, \tilde{l}_3)} . 
\end{equation*}
Hence we obtain Eq.(\ref{eq:HL0dL0dH}).
We will prove (ii) in the appendix.
\end{proof}

Let $f_1(x,E)$ and $f_2 (x,E)$ be a basis of solutions to the differential equation $(H^{(l_0,l_1,l_2,l_3)} -E) f(x)=0$.
Since the operator $H^{(l_0,l_1,l_2,l_3)}$ is doubly-periodic, the functions $f_1(x+2\omega _k,E)$ and $f_2 (x+2\omega _k,E)$ $(k=1,3)$ are also solutions to $(H^{(l_0,l_1,l_2,l_3)} -E) f(x)=0$, and they are written as linear combinations of $f_1(x,E)$ and $f_2 (x,E)$. Hence we have monodromy matrices
\begin{align}
& (f_1(x+2\omega _k,E) \; f_2(x+2\omega _k,E)) = (f_1(x,E) \; f_2(x,E)) 
\left(
\begin{array}{cc}
a_{11}^{(k)} & a_{12}^{(k)} \\
a_{21}^{(k)} & a_{22}^{(k)} 
\end{array}
\right) .
\label{eq:f1f2monod}
\end{align}
Now assume that $\alpha _0, \dots ,\alpha _3, d$ satisfy the assumption of Theorem \ref{thm:HL0123L0123H}. Set $\tilde{f}_i(x,E) =L_{\alpha _0, \alpha _1, \alpha _2, \alpha _3}f_i(x,E)$ $(i=1,2)$. It follows from Eq.(\ref{eq:HL0123L0123H}) that $\tilde{f}_1(x,E)$ and $\tilde{f}_2(x,E)$ are eigenfunctions of the operator $H^{(\alpha _0 +d,\alpha _1 +d,\alpha _2 +d,\alpha _3 +d)}$ with the eigenvalue $E$. If $E$ is not an eigenvalue of $H^{(l_0,l_1,l_2,l_3)} $ on the space $V_{\alpha _0, \alpha _1, \alpha _2, \alpha _3}$, then $\{ \tilde{f}_1(x,E)$, $\tilde{f}_2(x,E) \}$ is a basis of solutions to the differential equation $(H^{(\alpha _0 +d,\alpha _1 +d,\alpha _2 +d,\alpha _3 +d)} -E) f(x)=0$. It is shown in Proposition \ref{prop:L0123} that operator $L_{\alpha _0, \alpha _1, \alpha _2, \alpha _3}$ is doubly-periodic, and it follows from Eq.(\ref{eq:f1f2monod}) that
\begin{align*}
& (\tilde{f}_1(x+2\omega _k,E) \; \tilde{f}_2(x+2\omega _k,E)) = (\tilde{f}_1(x,E) \; \tilde{f}_2(x,E)) 
\left(
\begin{array}{cc}
a_{11}^{(k)} & a_{12}^{(k)} \\
a_{21}^{(k)} & a_{22}^{(k)} 
\end{array}
\right) . 
\end{align*}
Hence the monodromy structure of $H^{(l_0,l_1,l_2,l_3)}$ coincides with the one of $H^{(\alpha _0 +d,\alpha _1 +d,\alpha _2 +d,\alpha _3 +d)}$ for the case that $E$ is not an eigenvalue of $H^{(l_0,l_1,l_2,l_3)} $ on the space $V_{\alpha _0, \alpha _1, \alpha _2, \alpha _3}$. If $E$ is an eigenvalue of $H^{(l_0,l_1,l_2,l_3)} $ on the space $V_{\alpha _0, \alpha _1, \alpha _2, \alpha _3}$, it is also the eigenvalue of the operator $H^{(\alpha _0 +d,\alpha _1 +d,\alpha _2 +d,\alpha _3 +d)}$ on the space $V_{-\alpha _0-d, -\alpha _1-d, -\alpha _2-d, -\alpha _3-d}$, which follows from Proposition \ref{prop:HL0dL0dH} (ii).
Thus the operator $L_{\alpha _0, \alpha _1, \alpha _2,\alpha _3}$ defines an isomonodromic transformation from $H^{(l_0,l_1,l_2,l_3)}$ to $H^{(\alpha _0 +d,\alpha _1 +d,\alpha _2 +d,\alpha _3 +d)}$. Similarly  $L_{-\alpha _0-d, -\alpha _1-d, -\alpha _2-d, -\alpha _3-d}$ defines an isomonodromic transformation from $H^{(\alpha _0 +d,\alpha _1 +d,\alpha _2 +d,\alpha _3 +d)}$ to $H^{(l_0,l_1,l_2,l_3)}$.

On the multiplicity of eigenvalues of the operator $H^{(l_0,l_1,l_2,l_3)}$ on the space $V _{\alpha _0, \alpha _1, \alpha _2, \alpha _3}$, we have the following proposition:
\begin{prop} \label{prop:distrt} (c.f. \cite{TakR})
Let $\alpha_i$ be a number such that $\alpha_i= -l_i$ or $\alpha_i= l_i+1$ for each $i\in \{ 0,1,2,3\} $. Set $d=-\sum_{i=0}^3 \alpha_i /2$ and assume $d\in \Zint_{\geq 0}$ and $\alpha _i \neq \alpha _j$ for some $i,j \in \{0,1,2,3\}$. Then zeros of the characteristic polynomial of the operator $H^{(l_0,l_1,l_2,l_3)}$ on the space $V _{\alpha _0, \alpha _1, \alpha _2, \alpha _3}$ are distinct for generic periods $(2\omega_1, 2\omega _3)$.
\end{prop}
We prove this proposition in the appendix. Now we give a remark on the meaning of ``generic''. In the appendix, it is shown that there exists a periods $(2\omega _1, 2\omega _3)$ such that zeros of the characteristic polynomial are distinct. Hence the discriminant of the characteristic polynomial is not identically zero, and the set of periods such that the discriminant of the characteristic polynomial is not equal to zero is open-dense.

\section{Integer case} \label{sec:int}
We investigate quasi-solvability and generalized Darboux transformation for the case $l_i \in \Zint_{\geq 0}$ $(i=0,1,2,3)$.
Throughout this section, assume $l_i \in \Zint_{\geq 0}$ for $i=0,1,2,3$.

Let ${\mathcal F}$ be the space spanned by meromorphic doubly periodic functions up to signs, namely
\begin{align*}
& {\mathcal F}=\bigoplus _{\epsilon _1 , \epsilon _3 =\pm 1 } {\mathcal F} _{\epsilon _1 , \epsilon _3 }, \\
& {\mathcal F} _{\epsilon _1 , \epsilon _3 }=\{ f(x) \mbox{: meromorphic }| f(x+2\omega_1)= \epsilon _1 f(x), \; f(x+2\omega_3)= \epsilon _3 f(x) \} ,
\end{align*}
If $\alpha _i \in \Zint$ $(i=0,1,2,3)$ and $-\sum_{i=0}^3 \alpha _i/2 \in \Zint_{\geq 0}$, then $V_{\alpha _0, \alpha _1, \alpha _2, \alpha _3}$ is a subspace of ${\mathcal F} _{\epsilon _1 , \epsilon _3 }$ for suitable ${\mathcal F} _{\epsilon _1 , \epsilon _3 }$ $(\epsilon _1, \epsilon _3 \in \{ \pm 1 \})$, because $(\wp (x)-e_i)^{1/2}= \wp _i (x)$ $(i=1,2,3)$, where $\wp _i (x)$ is the co-$\wp $ function.
Invariant subspaces of ${\mathcal F}$ with respect to the operator $H^{(l_0,l_1,l_2,l_3)}$ is studied in \cite{Tak1} (see also \cite{GW,Tak2,Tak3}).

Let  $\alpha _i \in \{ -l_i , l_i+1 \}$ $(i=0,1,2,3)$,
\begin{align}
& U_{\alpha _0, \alpha _1, \alpha _2, \alpha _3}=
\left\{
\begin{array}{ll}
V_{\alpha _0, \alpha _1, \alpha _2, \alpha _3}, & \sum_{i=0}^3 \alpha _i/2 \in \Zint_{\leq 0} ;\\
V_{1-\alpha _0, 1-\alpha _1, 1-\alpha _2, 1-\alpha _3} ,& \sum_{i=0}^3 \alpha _i /2\in \Zint_{\geq 2} ;\\
\{ 0 \} ,& \mbox{otherwise},
\end{array}
\right. \label{eq:U}
\end{align}
and 
\begin{align}
& \tilde{L}_{\alpha _0, \alpha _1, \alpha _2, \alpha _3}=
\left\{
\begin{array}{ll}
L_{\alpha _0, \alpha _1, \alpha _2, \alpha _3}, & \sum_{i=0}^3 \alpha _i/2 \in \Zint_{\leq 0} ;\\
L_{1-\alpha _0, 1-\alpha _1, 1-\alpha _2, 1-\alpha _3} ,& \sum_{i=0}^3 \alpha _i /2\in \Zint_{\geq 2} ;\\
1 ,& \mbox{otherwise}.
\end{array}
\right. \label{eq:tL}
\end{align}

If $l_0 +l_1 +l_2 +l_3$ is even, then the operator $H^{(l_0,l_1,l_2,l_3)}$ (see Eq.(\ref{Ino})) preserves the spaces
\begin{equation}
U_{-l_0,-l_1,-l_2,-l_3}, \; U_{-l_0 ,-l_1,l_2+1 ,l_3+1}, \; U_{-l_0,l_1+1,-l_2,l_3+1}, \; U_{-l_0 ,l_1+1,l_2+1,-l_3}.
\label{eq:Ve0}
\end{equation}
Each space is contained in ${\mathcal F} _{\epsilon _1 , \epsilon _3 }$ for some $\epsilon _1 , \epsilon _3 \in \{\pm 1\}$, and the correspondence between the spaces and the signs of $(\epsilon _1 , \epsilon _3 )$ is one-to-one.
Let $V$ be the sum of these spaces.
Then $V$ is written as the direct sum of these spaces, i.e.
\begin{equation}
V=U_{-l_0,-l_1,-l_2,-l_3}\oplus U_{-l_0 ,-l_1,l_2+1 ,l_3+1}\oplus  U_{-l_0,l_1+1,-l_2,l_3+1}\oplus  U_{-l_0 ,l_1+1,l_2+1,-l_3}, \label{eq:Ve}
\end{equation}
and it is the maximal finite-dimensional invariant subspace in ${\mathcal F}$ with respect to the action of the operator $H^{(l_0,l_1,l_2,l_3)}$.
Let $k_i$ be the rearrangement of $l_i$ such that $k_0\geq k_1 \geq k_2 \geq k_3 (\geq 0)$. If $k_0+k_3\geq k_1+k_2$ (resp. $k_0+k_3< k_1+k_2$), then the dimension of the space $V$ is equal to $2k_0+1$ (resp. $k_0+k_1+k_2-k_3+1$).
Set 
\begin{align}
& g=
\left\{
\begin{array}{ll}
k_0, & k_0+k_3 \geq k_1+k_2 ;\\
(k_0+k_1+k_2-k_3)/2,& k_0+k_3 < k_1+k_2 .
\end{array}
\right. 
\label{eq:ge}
\end{align}
Then $g \in \Zint _{\geq 0}$ and $\dim V=2g+1$. We set
\begin{align}
& l_0^e= (-l_0+l_1+l_2+l_3)/2, \quad l_1^e= (l_0-l_1+l_2+l_3)/2, \label{def:lie} \\
& l_2^e= (l_0+l_1-l_2+l_3)/2, \quad l_3^e= (l_0+l_1+l_2-l_3)/2. \nonumber
\end{align}
Note that $l_0^e+ l_1^e +l_2^e+ l_3^e = l_0+l_1+l_2+l_3 \in 2 \Zint$.
It follows directly from Proposition \ref{eq:HL0123L0123H} that
\begin{align*}
& H^{(l_0^e,l_1^e,l_2^e,l_3^e)} \tilde{L}_{-l_0,-l_1,-l_2,-l_3}=\tilde{L}_{-l_0,-l_1,-l_2,-l_3} H^{(l_0,l_1,l_2,l_3)}, \\
& H^{(l_1^e,l_0^e,l_3^e,l_2^e)} \tilde{L}_{-l_0 ,-l_1,l_2+1 ,l_3+1}=\tilde{L}_{-l_0 ,-l_1,l_2+1 ,l_3+1} H^{(l_0,l_1,l_2,l_3)}, \nonumber \\
& H^{(l_2^e,l_3^e,l_0^e,l_1^e)} \tilde{L}_{-l_0,l_1+1,-l_2,l_3+1}=\tilde{L}_{-l_0,l_1+1,-l_2,l_3+1} H^{(l_0,l_1,l_2,l_3)}, \nonumber \\
& H^{(l_3^e,l_2^e,l_1^e,l_0^e)} \tilde{L}_{-l_0 ,l_1+1,l_2+1,-l_3}=\tilde{L}_{-l_0 ,l_1+1,l_2+1,-l_3} H^{(l_0,l_1,l_2,l_3)}. \nonumber 
\end{align*}
Thus it is shown that the operators which are linked by generalized Darboux transformations from $H^{(l_0,l_1,l_2,l_3)}$ are 
\begin{align}
& H^{(l_0^e,l_1^e,l_2^e,l_3^e)}, \; H^{(l_1^e,l_0^e,l_3^e,l_2^e)}, \; H^{(l_2^e,l_3^e,l_0^e,l_1^e)}, \; H^{(l_3^e,l_2^e,l_1^e,l_0^e)}, \label{eq:gDHe} \\
& H^{(l_0,l_1,l_2,l_3)}, \; H^{(l_1,l_0,l_3,l_2)}, \; H^{(l_2,l_3,l_0,l_1)}, \; H^{(l_3,l_2,l_1,l_0)}. \nonumber
\end{align}
As is discussed in section \ref{sec:Heun}, these eight operators are isomonodromic.
Note that the operators $H^{(l_1,l_0,l_3,l_2)}, \; H^{(l_2,l_3,l_0,l_1)}, \; H^{(l_3,l_2,l_1,l_0)}$ are obtained by the shift $x \rightarrow x+\omega _i$ $(i=1,2,3)$ from the operator $H^{(l_0,l_1,l_2,l_3)}$.

Assume that $l_0, \; l_1, \; l_2, \; l_3 \in \Zint$, $l_0 \geq l_1 \geq l_2 \geq l_3 \geq 0$ and $l_0 +l_1 +l_2 +l_3$ is even. 
Set $\tilde{l} _0 = l_3^e$, $\tilde{l} _1 = l_2^e$, $\tilde{l} _2 = l_1^e$ and $\tilde{l} _3 = \max ( l_0^e, -l_0^e-1)$. Then $\tilde{l} _0 \geq \tilde{l} _1 \geq \tilde{l} _2 \geq \tilde{l} _3 \geq 0$, and the operator $H^{(l_0,l_1,l_2,l_3)}$ is isomonodromic to $H^{(\tilde{l} _0, \tilde{l} _1, \tilde{l} _2, \tilde{l} _3)}$.
Note that, if $l_0+l_3 \neq l_1+l_2$ (resp. $l_0+l_3 = l_1+l_2$), then we have $(l_0,l_1,l_2,l_3)\neq (\tilde{l} _0, \tilde{l} _1, \tilde{l} _2, \tilde{l} _3)$ (resp. $(l_0,l_1,l_2,l_3) =(\tilde{l} _0, \tilde{l} _1, \tilde{l} _2, \tilde{l} _3)$).

We consider the case $l_0 +l_1 +l_2 +l_3$: odd.  If $l_0 +l_1 +l_2 +l_3$ is odd, then the operator $H^{(l_0,l_1,l_2,l_3)}$ preserves the spaces 
\begin{equation}
U_{-l_0,-l_1,-l_2,l_3+1}, \; U_{-l_0 ,-l_1,l_2+1,-l_3}, \; U_{-l_0 ,l_1+1,-l_2,-l_3}, \; U_{l_0+1,-l_1,-l_2,-l_3}.
\label{eq:Vo0}
\end{equation}
Each space is contained in ${\mathcal F} _{\epsilon _1 , \epsilon _3 }$ for some $\epsilon _1 , \epsilon _3 \in \{\pm 1\}$, and the correspondence between the spaces and the signs of $(\epsilon _1 , \epsilon _3 )$ is one-to-one.
Let $V$ be the sum of these spaces. Then $V$ is written as
\begin{equation}
V=U_{-l_0,-l_1,-l_2,l_3+1}\oplus U_{-l_0 ,-l_1,l_2+1,-l_3}\oplus U_{-l_0 ,l_1+1,-l_2,-l_3}\oplus U_{l_0+1,-l_1,-l_2,-l_3}, \label{eq:Vo}
\end{equation}
and it is the maximal finite-dimensional invariant subspace in ${\mathcal F}$ with respect to the action of the operator $H^{(l_0,l_1,l_2,l_3)}$.
Let $k_i$ be the rearrangement of $l_i$ such that $k_0\geq k_1 \geq k_2 \geq k_3 (\geq 0)$. If $k_0\geq k_1+k_2+k_3+1$ (resp. $k_0 < k_1+k_2+k_3+1$), then the dimension of the space $V$ is equal to $2k_0+1$ (resp. $k_0+k_1+k_2+k_3+2$).
Set 
\begin{align}
& g=
\left\{
\begin{array}{ll}
k_0, & k_0 \geq k_1+k_2 +k_3+1;\\
(k_0+k_1+k_2+k_3+1)/2,& k_0 < k_1+k_2+k_3+1 .
\end{array}
\right. 
\label{eq:go}
\end{align}
Then $g \in \Zint _{\geq 0}$ and $\dim V=2g+1$. We set 
\begin{align}
& l_0^o= (l_0+l_1+l_2+l_3+1)/2, \quad l_1^o= (l_0+l_1-l_2-l_3-1)/2, \label{def:lio}\\
& l_2^o= (l_0-l_1+l_2-l_3-1)/2, \quad l_3^o= (l_0-l_1-l_2+l_3-1)/2. \nonumber
\end{align}
It follows from Proposition \ref{eq:HL0123L0123H} that
\begin{align*}
& H^{(l_0^o,l_1^o,l_2^o,l_3^o)} \tilde{L}_{l_0+1,-l_1,-l_2,-l_3}=\tilde{L}_{l_0+1,-l_1,-l_2,-l_3} H^{(l_0,l_1,l_2,l_3)}, \\
& H^{(l_1^o,l_0^o,l_3^o,l_2^o)} \tilde{L}_{-l_0 ,l_1+1,-l_2 ,-l_3}=\tilde{L}_{-l_0 ,l_1+1,-l_2 ,-l_3} H^{(l_0,l_1,l_2,l_3)}, \nonumber \\
& H^{(l_2^o,l_3^o,l_0^o,l_1^o)} \tilde{L}_{-l_0,-l_1,l_2 +1,-l_3}=\tilde{L}_{-l_0,-l_1,l_2+1,-l_3} H^{(l_0,l_1,l_2,l_3)}, \nonumber \\
& H^{(l_3^o,l_2^o,l_1^o,l_0^o)} \tilde{L}_{-l_0 ,-l_1,-l_2,l_3+1}=\tilde{L}_{-l_0 ,-l_1,-l_2,l_3+1} H^{(l_0,l_1,l_2,l_3)} .\nonumber 
\end{align*}
It is shown that the operators which are linked by generalized Darboux transformations from $ H^{(l_0,l_1,l_2,l_3)}$ are 
\begin{align}
& H^{(l_0^o,l_1^o,l_2^o,l_3^o)}, \; H^{(l_1^o,l_0^o,l_3^o,l_2^o)}, \; H^{(l_2^o,l_3^o,l_0^o,l_1^o)}, \; H^{(l_3^o,l_2^o,l_1^o,l_0^o)}, \label{eq:gDHo} \\
& H^{(l_0,l_1,l_2,l_3)}, \; H^{(l_1,l_0,l_3,l_2)}, \; H^{(l_2,l_3,l_0,l_1)}, \; H^{(l_3,l_2,l_1,l_0)}, \nonumber
\end{align}
and these eight operators are isomonodromic.

Assume that $l_0, \; l_1, \; l_2, \; l_3 \in \Zint $, $l_0 \geq l_1 \geq l_2 \geq l_3 \geq 0$ and $l_0 +l_1 +l_2 +l_3$ is odd. Set $\tilde{l} _0 = l_0^o$, $\tilde{l} _1 = l_1^o$, $\tilde{l} _2 = l_2^o$ and $\tilde{l} _3 = \max ( l_3^o, -l_3^o-1)$. Then $\tilde{l} _0 \geq \tilde{l} _1 \geq \tilde{l} _2 \geq \tilde{l} _3 \geq 0$, and the operator $H^{(l_0,l_1,l_2,l_3)}$ is isomonodromic to $H^{(\tilde{l} _0, \tilde{l} _1, \tilde{l} _2, \tilde{l} _3)}$.
Note that, if $l_0 \neq l_1+l_2+l_3+1$ (resp. $l_0 = l_1+l_2+l_3+1$), then we have $(l_0,l_1,l_2,l_3)\neq (\tilde{l} _0, \tilde{l} _1, \tilde{l} _2, \tilde{l} _3)$ (resp. $(l_0,l_1,l_2,l_3) =(\tilde{l} _0, \tilde{l} _1, \tilde{l} _2, \tilde{l} _3)$).

Now we reproduce results in this section for the case of Lam\'e equation and associated Lam\'e equation.
If three (resp. two) of $l_i$ ($i=0,1,2,3$) are zero, then Eq.(\ref{Heq}) is called the Lam\'e equation (resp. the associated Lam\'e equation). For simplicity, we consider the case $l_2=l_3=0$, $l_0, l_1 \in \Zint$ and $l_0 \geq l_1 \geq 0$.

If $l_0 +l_1$ is even and $l_0>l_1$, then the operator $H^{(l_0,l_1,0,0)}$ is isomonodromic to $H^{((l_0+l_1)/2, (l_0+l_1)/2,(l_0-l_1)/2,(l_0-l_1)/2-1)}$ by the transformation $L_{-l_0,l_1+1,1,0}$. Especially, if $l_0$ is even and $l_1=0$ (the case of Lam\'e equation), then $H^{(l_0,0,0,0)}$ is isomonodromic to $H^{(l_0/2, l_0/2,l_0/2,l_0/2-1)}$. Note that, if $l_0 +l_1$ is even and $l_0=l_1$, then the operator $H^{(l_0,l_1,0,0)}$ is self-dual.

If $l_0 +l_1$ is odd and $l_0>l_1-1$, then the operator $H^{(l_0,l_1,0,0)}$ is isomonodromic to $H^{((l_0+l_1+1)/2, (l_0+l_1-1)/2,(l_0-l_1-1)/2,(l_0-l_1-1)/2)}$ by the transformation $L_{-l_0,l_1+1,1,1}$. Especially, if $l_0$ is odd and $l_1=0$ (the case of Lam\'e equation), then $H^{(l_0,0,0,0)}$ is isomonodromic to $H^{((l_0+1)/2, (l_0-1)/2,(l_0-1)/2,(l_0-1)/2)}$. Note that, if $l_0 +l_1$ is odd and $l_0=l_1-1$, then the operator $H^{(l_0,l_1,0,0)}$ is self-dual.

Thus we confirmed the conjecture of Khare and Sukhatme \cite{KS} on isospectral  partner of Lam\'e equation and associated Lam\'e equation.

\section{Half-integer case} \label{sec:hint}
We investigate quasi-solvability and generalized Darboux transformation for the case $l_i + 1/2 \in \Zint  $ $(i=0,1,2,3)$. Write $l_i=n_i -1/2$ and assume $n_i \in \Zint_{\geq 0}$ for $i=0,1,2,3$.

Candidates for the space related to quasi-solvability written as Eq.(\ref{eq:Vqsol}) are described as $V_{\tilde{n}_0+\frac{1}{2},\tilde{n}_1+\frac{1}{2},\tilde{n}_2+\frac{1}{2},\tilde{n}_3+\frac{1}{2}}$, where $\tilde{n}_i \in \{ \pm n_i \}$ $(i=0,1,2,3)$.
For the existence of non-zero space $V_{\tilde{n}_0+\frac{1}{2},\tilde{n}_1+\frac{1}{2},\tilde{n}_2+\frac{1}{2},\tilde{n}_3+\frac{1}{2}}$, the condition $\sum _{i=0}^3 \tilde{n}_i \in 2 \Zint$ is necessary.
In other words, if $\sum _{i=0}^3 n_i$ is odd, then the space related to quasi-solvability written as Eq.(\ref{eq:Vqsol}) does not exist.

On the rest of this section, we assume that $\sum _{i=0}^3 n_i$ is even. Then the operator $H^{(n_0-\frac{1}{2},n_1-\frac{1}{2},n_2-\frac{1}{2},n_3-\frac{1}{2})}$ preserves the spaces
\begin{align*}
& U_{-n_0+\frac{1}{2},-n_1+\frac{1}{2},-n_2+\frac{1}{2},-n_3+\frac{1}{2}}, \; U_{-n_0+\frac{1}{2},-n_1+\frac{1}{2},n_2+\frac{1}{2},n_3+\frac{1}{2}}, \\
&  U_{-n_0+\frac{1}{2},n_1+\frac{1}{2},-n_2+\frac{1}{2},n_3+\frac{1}{2}}, \; U_{-n_0+\frac{1}{2},n_1+\frac{1}{2},n_2+\frac{1}{2},-n_3+\frac{1}{2}},  \nonumber \\
&  U_{-n_0+\frac{1}{2},-n_1+\frac{1}{2},-n_2+\frac{1}{2},n_3+\frac{1}{2}}, \; U_{-n_0+\frac{1}{2},-n_1+\frac{1}{2},n_2+\frac{1}{2},-n_3+\frac{1}{2}}, \nonumber \\
& U_{-n_0+\frac{1}{2},n_1+\frac{1}{2},-n_2+\frac{1}{2},-n_3+\frac{1}{2}}, \; U_{n_0+\frac{1}{2},-n_1+\frac{1}{2},-n_2+\frac{1}{2},-n_3+\frac{1}{2}}, \nonumber
\end{align*}
where $U_{\alpha _0, \alpha _1, \alpha _2, \alpha _3}$ is defined in Eq.(\ref{eq:U}). Unlike the case $l_0 ,l_1 ,l_2 ,l_3 \in \Zint$, these eight spaces are not disjoint. In fact, these eight spaces are subspaces of the space $U_{-n_0+\frac{1}{2},-n_1+\frac{1}{2},-n_2+\frac{1}{2},-n_3+\frac{1}{2}}$. Set
\begin{align*}
& n_0^{(1)}= (n_0+n_1+n_2+n_3)/2, \quad n_1^{(1)}= (n_0+n_1-n_2-n_3)/2, \\
& n_2^{(1)}= (n_0-n_1+n_2-n_3)/2, \quad n_3^{(1)}= (n_0-n_1-n_2+n_3)/2, \nonumber \\
& n_0^{(2)}= (-n_0+n_1+n_2+n_3)/2, \quad n_1^{(2)}= (n_0-n_1+n_2+n_3)/2, \nonumber \\
& n_2^{(2)}= (n_0+n_1-n_2+n_3)/2, \quad n_3^{(2)}= (n_0+n_1+n_2-n_3)/2. \nonumber
\end{align*}
It follows from Proposition \ref{eq:HL0123L0123H} that
\begin{align*}
& H^{(n_0^{(1)}-\frac{1}{2},n_1^{(1)}-\frac{1}{2},n_2^{(1)}-\frac{1}{2},n_3^{(1)}-\frac{1}{2})} \tilde{L}_{n_0+\frac{1}{2},-n_1+\frac{1}{2},-n_2+\frac{1}{2},-n_3+\frac{1}{2}} \\
& \; =\tilde{L}_{n_0+\frac{1}{2},-n_1+\frac{1}{2},-n_2+\frac{1}{2},-n_3+\frac{1}{2}} H^{(n_0-\frac{1}{2},n_1-\frac{1}{2},n_2-\frac{1}{2},n_3-\frac{1}{2})}, \nonumber \\
& H^{(n_1^{(1)}-\frac{1}{2},n_0^{(1)}-\frac{1}{2},n_3^{(1)}-\frac{1}{2},n_2^{(1)}-\frac{1}{2})} \tilde{L}_{-n_0+\frac{1}{2},n_1+\frac{1}{2},-n_2+\frac{1}{2},-n_3+\frac{1}{2}} \nonumber \\
& \; =\tilde{L}_{-n_0+\frac{1}{2},n_1+\frac{1}{2},-n_2+\frac{1}{2},-n_3+\frac{1}{2}} H^{(n_0-\frac{1}{2},n_1-\frac{1}{2},n_2-\frac{1}{2},n_3-\frac{1}{2})}, \nonumber \\
& H^{(n_2^{(1)}-\frac{1}{2},n_3^{(1)}-\frac{1}{2},n_0^{(1)}-\frac{1}{2},n_1^{(1)}-\frac{1}{2})} \tilde{L}_{-n_0+\frac{1}{2},-n_1+\frac{1}{2},n_2+\frac{1}{2},-n_3+\frac{1}{2}} \nonumber \\
& \; =\tilde{L}_{-n_0+\frac{1}{2},-n_1+\frac{1}{2},n_2+\frac{1}{2},-n_3+\frac{1}{2}} H^{(n_0-\frac{1}{2},n_1-\frac{1}{2},n_2-\frac{1}{2},n_3-\frac{1}{2})}, \nonumber \\
& H^{(n_3^{(1)}-\frac{1}{2},n_2^{(1)}-\frac{1}{2},n_1^{(1)}-\frac{1}{2},n_0^{(1)}-\frac{1}{2})} \tilde{L}_{-n_0+\frac{1}{2},-n_1+\frac{1}{2},-n_2+\frac{1}{2},n_3+\frac{1}{2}} \nonumber \\
& \; =\tilde{L}_{-n_0+\frac{1}{2},-n_1+\frac{1}{2},-n_2+\frac{1}{2},n_3+\frac{1}{2}} H^{(n_0-\frac{1}{2},n_1-\frac{1}{2},n_2-\frac{1}{2},n_3-\frac{1}{2})}, \nonumber 
\end{align*}
where $\tilde{L}_{\alpha _0, \alpha _1, \alpha _2, \alpha _3}$ is defined in Eq.(\ref{eq:tL}). We also obtain that
\begin{align*}
& H^{(n_0^{(2)}-\frac{1}{2},n_1^{(2)}-\frac{1}{2},n_2^{(2)}-\frac{1}{2},n_3^{(2)}-\frac{1}{2})} \tilde{L}_{-n_0+\frac{1}{2},-n_1+\frac{1}{2},-n_2+\frac{1}{2},-n_3+\frac{1}{2}} \\
& \; =\tilde{L}_{-n_0+\frac{1}{2},-n_1+\frac{1}{2},-n_2+\frac{1}{2},-n_3+\frac{1}{2}} H^{(n_0-\frac{1}{2},n_1-\frac{1}{2},n_2-\frac{1}{2},n_3-\frac{1}{2})}, \nonumber \\
& H^{(n_1^{(2)}-\frac{1}{2},n_0^{(2)}-\frac{1}{2},n_3^{(2)}-\frac{1}{2},n_2^{(2)}-\frac{1}{2})} \tilde{L}_{-n_0+\frac{1}{2},-n_1+\frac{1}{2},n_2+\frac{1}{2},n_3+\frac{1}{2}} \nonumber \\
& \; =\tilde{L}_{-n_0+\frac{1}{2},-n_1+\frac{1}{2},n_2+\frac{1}{2},n_3+\frac{1}{2}} H^{(n_0-\frac{1}{2},n_1-\frac{1}{2},n_2-\frac{1}{2},n_3-\frac{1}{2})}, \nonumber \\
& H^{(n_2^{(2)}-\frac{1}{2},n_3^{(2)}-\frac{1}{2},n_0^{(2)}-\frac{1}{2},n_1^{(2)}-\frac{1}{2})} \tilde{L}_{-n_0+\frac{1}{2},n_1+\frac{1}{2},-n_2+\frac{1}{2},n_3+\frac{1}{2}} \nonumber \\
& \; =\tilde{L}_{-n_0+\frac{1}{2},n_1+\frac{1}{2},-n_2+\frac{1}{2},n_3+\frac{1}{2}} H^{(n_0-\frac{1}{2},n_1-\frac{1}{2},n_2-\frac{1}{2},n_3-\frac{1}{2})}, \nonumber \\
& H^{(n_3^{(2)}-\frac{1}{2},n_2^{(2)}-\frac{1}{2},n_1^{(2)}-\frac{1}{2},n_0^{(2)}-\frac{1}{2})} \tilde{L}_{-n_0+\frac{1}{2},n_1+\frac{1}{2},n_2+\frac{1}{2},-n_3+\frac{1}{2}} \nonumber \\
& \; =\tilde{L}_{-n_0+\frac{1}{2},n_1+\frac{1}{2},n_2+\frac{1}{2},-n_3+\frac{1}{2}} H^{(n_0-\frac{1}{2},n_1-\frac{1}{2},n_2-\frac{1}{2},n_3-\frac{1}{2})}. \nonumber 
\end{align*}
Thus it is shown that the operators which are linked by generalized Darboux transformations from $H^{(n_0-\frac{1}{2},n_1-\frac{1}{2},n_2-\frac{1}{2},n_3-\frac{1}{2})}$ are 
\begin{align*}
& H^{(n_0^{(1)}-\frac{1}{2},n_1^{(1)}-\frac{1}{2},n_2^{(1)}-\frac{1}{2},n_3^{(1)}-\frac{1}{2})}, \; H^{(n_1^{(1)}-\frac{1}{2},n_0^{(1)}-\frac{1}{2},n_3^{(1)}-\frac{1}{2},n_2^{(1)}-\frac{1}{2})}, \nonumber \\
& H^{(n_2^{(1)}-\frac{1}{2},n_3^{(1)}-\frac{1}{2},n_1^{(1)}-\frac{1}{2},n_0^{(1)}-\frac{1}{2})}, \; H^{(n_3^{(1)}-\frac{1}{2},n_2^{(1)}-\frac{1}{2},n_1^{(1)}-\frac{1}{2},n_0^{(1)}-\frac{1}{2})}, \nonumber \\
& H^{(n_0^{(2)}-\frac{1}{2},n_1^{(2)}-\frac{1}{2},n_2^{(2)}-\frac{1}{2},n_3^{(2)}-\frac{1}{2})}, \; H^{(n_1^{(2)}-\frac{1}{2},n_0^{(2)}-\frac{1}{2},n_3^{(2)}-\frac{1}{2},n_2^{(2)}-\frac{1}{2})}, \nonumber \\
& H^{(n_2^{(2)}-\frac{1}{2},n_3^{(2)}-\frac{1}{2},n_1^{(2)}-\frac{1}{2},n_0^{(2)}-\frac{1}{2})}, \; H^{(n_3^{(2)}-\frac{1}{2},n_2^{(2)}-\frac{1}{2},n_1^{(2)}-\frac{1}{2},n_0^{(2)}-\frac{1}{2})}, \nonumber \\
& H^{(n_0-\frac{1}{2},n_1-\frac{1}{2},n_2-\frac{1}{2},n_3-\frac{1}{2})}, \; H^{(n_1-\frac{1}{2},n_0-\frac{1}{2},n_3-\frac{1}{2},n_2-\frac{1}{2})}, \\
& H^{(n_2-\frac{1}{2},n_3-\frac{1}{2},n_0-\frac{1}{2},n_1-\frac{1}{2})}, \; H^{(n_3-\frac{1}{2},n_2-\frac{1}{2},n_1-\frac{1}{2},n_0-\frac{1}{2})}. \nonumber 
\end{align*}
Now we investigate the case $n_0, \; n_1, \; n_2, \; n_3 \in \Zint _{\geq 0}$, $n_0+ n_1+ n_2+ n_3 \in 2\Zint _{\geq 0}$ and $n_0 \geq n_1 \geq n_2 \geq n_3 $. Recall that $l_i=n_i-1/2$ ($i=0,1,2,3$).
We divide into three cases, the case $n_0 \geq n_1 +n_2+n_3$, the case $n_1+n_2-n_3 \leq n_0 < n_1 +n_2+n_3$ and the case $n_0< n_1 +n_2 -n_3$.
if $\alpha _0 +\alpha _1 +\alpha _2 +\alpha _3=2$, then we set $V _{\alpha _0 ,\alpha _1 ,\alpha _2 ,\alpha _3}=\{0\}$.

If $n_0 \geq n_1 +n_2+n_3$, then the operator $H^{(n_0-\frac{1}{2},n_1-\frac{1}{2},n_2-\frac{1}{2},n_3-\frac{1}{2})}$ preserves the spaces
\begin{center}
\begin{picture}(430,80)(0,0)
\put(0,55){$V_{\frac{1}{2}-n_0,\frac{1}{2}-n_1,\frac{1}{2}-n_2,\frac{1}{2}-n_3}$}
\put(85,70){$V_{\frac{1}{2}-n_0,\frac{1}{2}-n_1,\frac{1}{2}-n_2,\frac{1}{2}+n_3}$}
\put(85,40){$V_{\frac{1}{2}-n_0,\frac{1}{2}-n_1,\frac{1}{2}+n_2,\frac{1}{2}-n_3}$}
\put(85,10){$V_{\frac{1}{2}-n_0,\frac{1}{2}+n_1,\frac{1}{2}-n_2,\frac{1}{2}-n_3}$}
\put(225,70){$V_{\frac{1}{2}-n_0,\frac{1}{2}-n_1,\frac{1}{2}+n_2,\frac{1}{2}+n_3}$}
\put(225,40){$V_{\frac{1}{2}-n_0,\frac{1}{2}+n_1,\frac{1}{2}-n_2,\frac{1}{2}+n_3}$}
\put(225,10){$V_{\frac{1}{2}-n_0,\frac{1}{2}+n_1,\frac{1}{2}+n_2,\frac{1}{2}-n_3}$}
\put(325,25){$V_{\frac{1}{2}-n_0,\frac{1}{2}+n_1,\frac{1}{2}+n_2,\frac{1}{2}+n_3}$.}
\put(55,65){\line (3,1){15}}
\put(55,40){\line (1,-1){25}}
\put(55,45){\line (3,-1){15}}
\put(190,70){\line (1,0){30}}
\put(190,10){\line (1,0){30}}
\put(193,43){\line (1,1){24}}
\put(193,13){\line (1,1){24}}
\put(193,67){\line (1,-1){24}}
\put(193,37){\line (1,-1){24}}
\put(330,65){\line (1,-1){25}}
\put(330,40){\line (3,-1){15}}
\put(330,10){\line (3,1){15}}
\end{picture}
\end{center}
Here $V _{\alpha _0 ,\alpha _1 ,\alpha _2 ,\alpha _3}$---$V _{\tilde{\alpha }_0 ,\tilde{\alpha } _1 ,\tilde{\alpha } _2 ,\tilde{\alpha } _3}$ means that $V _{\tilde{\alpha }_0 ,\tilde{\alpha } _1 ,\tilde{\alpha } _2 ,\tilde{\alpha } _3}$ is a subspace of $V _{\alpha _0 ,\alpha _1 ,\alpha _2 ,\alpha _3}$.
On this case, the operator $H^{(n_0-\frac{1}{2},n_1-\frac{1}{2},n_2-\frac{1}{2},n_3-\frac{1}{2})}$ is isomonodromic to $H^{(n_0^{(1)}-\frac{1}{2},n_1^{(1)}-\frac{1}{2},n_2^{(1)}-\frac{1}{2},n_3^{(1)}-\frac{1}{2})}$ and $H^{(n_3^{(2)}-\frac{1}{2},n_2^{(2)}-\frac{1}{2},n_1^{(2)}-\frac{1}{2},-n_0^{(2)}-\frac{1}{2})}$ by the transformations $L_{-n_0+\frac{1}{2},n_1+\frac{1}{2},n_2+\frac{1}{2},n_3+\frac{1}{2}}$ and $L_{-n_0+\frac{1}{2},n_1+\frac{1}{2},n_2+\frac{1}{2},-n_3+\frac{1}{2}}$.

If $n_1+n_2-n_3 \leq n_0 < n_1 +n_2+n_3$, then the operator $H^{(n_0-\frac{1}{2},n_1-\frac{1}{2},n_2-\frac{1}{2},n_3-\frac{1}{2})}$ preserves the spaces
\begin{center}
\begin{picture}(430,80)(0,0)
\put(0,40){$V_{\frac{1}{2}-n_0,\frac{1}{2}-n_1,\frac{1}{2}-n_2,\frac{1}{2}-n_3}$}
\put(150,70){$V_{\frac{1}{2}-n_0,\frac{1}{2}-n_1,\frac{1}{2}-n_2,\frac{1}{2}+n_3}$}
\put(150,50){$V_{\frac{1}{2}-n_0,\frac{1}{2}-n_1,\frac{1}{2}+n_2,\frac{1}{2}-n_3}$}
\put(150,30){$V_{\frac{1}{2}-n_0,\frac{1}{2}+n_1,\frac{1}{2}-n_2,\frac{1}{2}-n_3}$}
\put(150,10){$V_{\frac{1}{2}+n_0,\frac{1}{2}-n_1,\frac{1}{2}-n_2,\frac{1}{2}-n_3}$}
\put(300,70){$V_{\frac{1}{2}-n_0,\frac{1}{2}-n_1,\frac{1}{2}+n_2,\frac{1}{2}+n_3}$}
\put(300,50){$V_{\frac{1}{2}-n_0,\frac{1}{2}+n_1,\frac{1}{2}-n_2,\frac{1}{2}+n_3}$}
\put(300,30){$V_{\frac{1}{2}-n_0,\frac{1}{2}+n_1,\frac{1}{2}+n_2,\frac{1}{2}-n_3}$,}
\put(115,46){\line (1,1){25}}
\put(115,42){\line (3,1){25}}
\put(115,38){\line (3,-1){25}}
\put(115,34){\line (1,-1){25}}
\put(260,75){\line (1,0){35}}
\put(260,25){\line (1,0){35}}
\put(260,70){\line (2,-1){35}}
\put(260,48){\line (2,-1){35}}
\put(260,52){\line (2,1){35}}
\put(260,30){\line (2,1){35}}
\end{picture}
\end{center}
and it is isomonodromic to $H^{(n_0^{(1)}-\frac{1}{2},n_1^{(1)}-\frac{1}{2},n_2^{(1)}-\frac{1}{2},n_3^{(1)}-\frac{1}{2})}$ and $H^{(n_3^{(2)}-\frac{1}{2},n_2^{(2)}-\frac{1}{2},n_1^{(2)}-\frac{1}{2},n_0^{(2)}-\frac{1}{2})}$ by the transformations $L_{n_0+\frac{1}{2},-n_1+\frac{1}{2},-n_2+\frac{1}{2},-n_3+\frac{1}{2}}$ and $L_{-n_0+\frac{1}{2},n_1+\frac{1}{2},n_2+\frac{1}{2},-n_3+\frac{1}{2}}$.

If $n_0< n_1 +n_2 -n_3$, then the operator $H^{(n_0-\frac{1}{2},n_1-\frac{1}{2},n_2-\frac{1}{2},n_3-\frac{1}{2})}$ preserves the spaces
\begin{center}
\begin{picture}(430,80)(0,0)
\put(0,40){$V_{\frac{1}{2}-n_0,\frac{1}{2}-n_1,\frac{1}{2}-n_2,\frac{1}{2}-n_3}$}
\put(150,70){$V_{\frac{1}{2}-n_0,\frac{1}{2}-n_1,\frac{1}{2}-n_2,\frac{1}{2}+n_3}$}
\put(150,50){$V_{\frac{1}{2}-n_0,\frac{1}{2}-n_1,\frac{1}{2}+n_2,\frac{1}{2}-n_3}$}
\put(150,30){$V_{\frac{1}{2}-n_0,\frac{1}{2}+n_1,\frac{1}{2}-n_2,\frac{1}{2}-n_3}$}
\put(150,10){$V_{\frac{1}{2}+n_0,\frac{1}{2}-n_1,\frac{1}{2}-n_2,\frac{1}{2}-n_3}$}
\put(300,60){$V_{\frac{1}{2}-n_0,\frac{1}{2}-n_1,\frac{1}{2}+n_2,\frac{1}{2}+n_3}$}
\put(300,40){$V_{\frac{1}{2}-n_0,\frac{1}{2}+n_1,\frac{1}{2}-n_2,\frac{1}{2}+n_3}$}
\put(300,20){$V_{\frac{1}{2}+n_0,\frac{1}{2}-n_1,\frac{1}{2}-n_2,\frac{1}{2}+n_3}$,}
\put(115,46){\line (1,1){25}}
\put(115,42){\line (3,1){25}}
\put(115,38){\line (3,-1){25}}
\put(115,34){\line (1,-1){25}}
\put(260,74){\line (4,-1){35}}
\put(260,71){\line (3,-2){35}}
\put(260,52){\line (4,1){35}}
\put(260,68){\line (4,-5){35}}
\put(260,32){\line (4,1){35}}
\put(260,12){\line (4,1){35}}
\end{picture}
\end{center}
and it is isomonodromic to $H^{(n_0^{(1)}-\frac{1}{2},n_1^{(1)}-\frac{1}{2},n_2^{(1)}-\frac{1}{2},-n_3^{(1)}-\frac{1}{2})}$ and $H^{(n_3^{(2)}-\frac{1}{2},n_2^{(2)}-\frac{1}{2},n_1^{(2)}-\frac{1}{2},n_0^{(2)}-\frac{1}{2})}$ by the transformations $L_{n_0+\frac{1}{2},-n_1+\frac{1}{2},-n_2+\frac{1}{2},-n_3+\frac{1}{2}}$ and $L_{n_0+\frac{1}{2},-n_1+\frac{1}{2},-n_2+\frac{1}{2},n_3+\frac{1}{2}}$.

\section{Finite-gap potential} \label{sec:fin}

In this section, we construct an odd-order differential operator $\tilde{A}$ (see below) which commutes with the operator $H^{(l_0,l_1,l_2,l_3)}$ by composing generalized Darboux transformations for the case $l_0, l_1, l_2, l_3 \in \Zint _{\geq 0}$. We also show that the commuting operator coincides with the one constructed in \cite{Tak3}.
We use notations in section \ref{sec:int}.
\begin{prop}
If $l_0+l_1+l_2+l_3$ is even, then we set
\begin{equation}
\tilde{A}= \tilde{L}_{-l_3^e,l_2^e+1,l_1^e+1,-l_0^e} \tilde{L}_{-l_1,l_0+1,-l_3,l_2+1} \tilde{L}_{-l_0^e,-l_1^e,l_2^e+1,l_3^e+1} \tilde{L}_{-l_0,-l_1,-l_2,-l_3} . \label{eq:Ate}
\end{equation}
If $l_0+l_1+l_2+l_3$ is odd, then we set
\begin{equation}
\tilde{A}= \tilde{L}_{l_2^o+1,-l_3^o,-l_0^o,-l_1^o} \tilde{L}_{-l_1,-l_0,l_3+1,-l_2} \tilde{L}_{-l_0^o,l_1^o+1,-l_2^o,-l_3^o} \tilde{L}_{l_0+1,-l_1,-l_2,-l_3} . \label{eq:Ato}
\end{equation}
The operator $\tilde {A}$ commutes with $H^{(l_0,l_1,l_2,l_3)}$, i.e.,
\begin{equation}
\tilde {A}H^{(l_0,l_1,l_2, l_3)}=H^{(l_0,l_1,l_2, l_3)}\tilde {A}. \label{eq:Acomme}
\end{equation}
\end{prop}
\begin{proof}
If $l_0+l_1+l_2+l_3$ is even, then Eq.(\ref{eq:Acomme}) is shown by applying the following relations:
\begin{align}
& H^{(l_0^e,l_1^e,l_2^e,l_3^e)} \tilde{L}_{-l_0,-l_1,-l_2,-l_3}=\tilde{L}_{-l_0,-l_1,-l_2,-l_3} H^{(l_0,l_1,l_2,l_3)},  \label{eq:HLLHe} \\
& H^{(l_1,l_0,l_3,l_2)} \tilde{L}_{-l_0^e,-l_1^e,l_2^e+1,l_3^e+1}=\tilde{L}_{-l_0^e,-l_1^e,l_2^e+1,l_3^e+1} H^{(l_0^e,l_1^e,l_2^e,l_3^e)}, \nonumber \\
& H^{(l_3^e,l_2^e,l_1^e,l_0^e)}\tilde{L}_{-l_1,l_0+1,-l_3,l_2+1} =\tilde{L}_{-l_1,l_0+1,-l_3,l_2+1} H^{(l_1,l_0,l_3,l_2)}, \nonumber \\
&  H^{(l_0,l_1,l_2,l_3)} \tilde{L}_{-l_3^e,l_2^e+1,l_1^e+1,-l_0^e}= \tilde{L}_{-l_3^e,l_2^e+1,l_1^e+1,-l_0^e}H^{(l_3^e,l_2^e,l_1^e,l_0^e)}. \nonumber 
\end{align}
If $l_0+l_1+l_2+l_3$ is odd, then Eq.(\ref{eq:Acomme}) is shown by
\begin{align}
& H^{(l_0^o,l_1^o,l_2^o,l_3^o)} \tilde{L}_{l_0+1,-l_1,-l_2,-l_3} = \tilde{L}_{l_0+1,-l_1,-l_2,-l_3} H^{(l_0,l_1,l_2,l_3)}, \\
& H^{(l_1,l_0,l_3,l_2)} \tilde{L}_{-l_0^o,l_1^o+1,-l_2^o,-l_3^o}=\tilde{L}_{-l_0^o,l_1^o+1,-l_2^o,-l_3^o} H^{(l_0^o,l_1^o,l_2^o,l_3^o)}, \nonumber \\
& H^{(l_2^o,l_3^o,l_0^o,l_1^o)} \tilde{L}_{-l_1,-l_0,l_3+1,-l_2} = \tilde{L}_{-l_1,-l_0,l_3+1,-l_2} H^{(l_1,l_0,l_3,l_2)}, \nonumber \\
&  H^{(l_0,l_1,l_2,l_3)} \tilde{L}_{l_2^o+1,-l_3^o,-l_0^o,-l_1^o}= \tilde{L}_{l_2^o+1,-l_3^o,-l_0^o,-l_1^o}H^{(l_2^o,l_3^o,l_0^o,l_1^o)}. \nonumber 
\end{align}
\end{proof}
Note that the order of the differential operator $\tilde{A}$ is equal to $2g+1$, where $g$ is defined in Eq.(\ref{eq:ge}) or Eq.(\ref{eq:go}).
Recall that the space $V$ is written as Eq.(\ref{eq:Ve}) or Eq.(\ref{eq:Vo}), $\dim V= 2g+1$ and it is the maximal finite-dimensional invariant subspace in ${\mathcal F}$ with respect to the action of the operator $H^{(l_0,l_1,l_2,l_3)}$.
\begin{prop} \label{prop:tAV}
The operator $\tilde{A}$ is the monic differential operator of minimum order which vanishes all elements in $V$.
\end{prop}
\begin{proof}
We show this proposition for the case that $l_0+l_1+l_2+l_3$ is even. For the case that $l_0+l_1+l_2+l_3$ is odd, it is shown similarly.

Since $V$ is written as Eq.(\ref{eq:Ve}) and $\deg \tilde{A} = \dim V$, it is enough to show that $\tilde{A} \phi (x) =0$ for $\phi (x) \in U_{-l_0,-l_1,-l_2,-l_3}$, $U_{-l_0 ,-l_1,l_2+1 ,l_3+1}$, $U_{-l_0,l_1+1,-l_2,l_3+1}$, $U_{-l_0 ,l_1+1,l_2+1,-l_3}$.
The operator $\tilde {A}$ is written as Eq.(\ref{eq:Ate}), and the operator $ \tilde{L}_{-l_0,-l_1,-l_2,-l_3} $ vanishes any element in the space $U_{-l_0,-l_1,-l_2,-l_3} $. Hence we have $\tilde{A} \phi (x) =0$ for $\phi (x) \in U_{-l_0,-l_1,-l_2,-l_3}$.
Let $\langle f_1 (x), \dots , f_{n} (x) \rangle $ be a basis of the space $U_{-l_0 ,-l_1,l_2+1 ,l_3+1}$. Since the operator $H^{(l_0,l_1,l_2, l_3)}$ preserve the space $U_{-l_0 ,-l_1,l_2+1 ,l_3+1}$, the function $H^{(l_0,l_1,l_2, l_3)} f_j (x)$ is written as $\sum _{i=1}^n a_{i,j} f _i(x)$ for some constants $a_{i,j} $.
It follows from Eq.(\ref{eq:HLLHe}) that 
\begin{align}
 H^{(l_0^e,l_1^e,l_2^e,l_3^e)} \tilde{L}_{-l_0,-l_1,-l_2,-l_3} f_j (x) & = \tilde{L}_{-l_0,-l_1,-l_2,-l_3} H^{(l_0,l_1,l_2,l_3)} f_j(x) \label{eq:Hl0eLt} \\
& =  \sum _{i=1}^n a_{i,j} \tilde{L}_{-l_0,-l_1,-l_2,-l_3} f_j (x). \nonumber 
\end{align}
Set $\tilde{U}_{-l_0^e,-l_1^e,l_2^e+1,l_3^e+1}= \tilde{L}_{-l_0,-l_1,-l_2,-l_3} U_{-l_0 ,-l_1,l_2+1 ,l_3+1}$. Then it follows from Eq.(\ref{eq:Hl0eLt}) that the space $\tilde{U}_{-l_0^e,-l_1^e,l_2^e+1,l_3^e+1}$ is invariant under the action of $H^{(l_0^e,l_1^e,l_2^e,l_3^e)}$.
Let $(\epsilon _1 , \epsilon _3 )$ $(\epsilon _1 , \epsilon _3 \in \{ \pm 1 \})$ be the numbers such that $U_{-l_0 ,-l_1,l_2+1 ,l_3+1} \subset {\mathcal F} _{\epsilon _1 , \epsilon _3 }$. Then we have $\tilde{U}_{-l_0^e,-l_1^e,l_2^e+1,l_3^e+1} \subset {\mathcal F} _{\epsilon _1 , \epsilon _3 }$ and $U_{-l_0^e,-l_1^e,l_2^e+1,l_3^e+1} \subset {\mathcal F} _{\epsilon _1 , \epsilon _3 }$.
As is shown in \cite[Theorem 3.1]{Tak1}, the space $U_{-l_0^e,-l_1^e,l_2^e+1,l_3^e+1}$ is the maximum subspace of ${\mathcal F} _{\epsilon _1 , \epsilon _3 }$ which is invariant under the action of $H^{(l_0^e,l_1^e,l_2^e,l_3^e)}$. Thus we have 
\begin{equation*}
\tilde{U}_{-l_0^e,-l_1^e,l_2^e+1,l_3^e+1} = \tilde{L}_{-l_0,-l_1,-l_2,-l_3} U_{-l_0 ,-l_1,l_2+1 ,l_3+1} \subset U_{-l_0^e,-l_1^e,l_2^e+1,l_3^e+1}.
\end{equation*}
Similarly it is shown that 
\begin{align*}
& \tilde{L}_{-l_0^e,-l_1^e,l_2^e+1,l_3^e+1} \tilde{L}_{-l_0,-l_1,-l_2,-l_3} U_{-l_0,l_1+1,-l_2,l_3+1} \subset U_{-l_1,l_0+1,-l_3,l_2+1} ,\\
& \tilde{L}_{-l_1,l_0+1,-l_3,l_2+1} \tilde{L}_{-l_0^e,-l_1^e,l_2^e+1,l_3^e+1} \tilde{L}_{-l_0,-l_1,-l_2,-l_3} U_{-l_0 ,l_1+1,l_2+1,-l_3}\subset U_{-l_3^e,l_2^e+1,l_1^e+1,-l_0^e}. \nonumber
\end{align*}
Since the operator $\tilde{L}_{-l_0^e,-l_1^e,l_2^e+1,l_3^e+1}$ (resp. $\tilde{L}_{-l_1,l_0+1,-l_3,l_2+1}$, $\tilde{L}_{-l_3^e,l_2^e+1,l_1^e+1,-l_0^e}$) vanishes any element in the space $U_{-l_0^e,-l_1^e,l_2^e+1,l_3^e+1}$ (resp. $U_{-l_1,l_0+1,-l_3,l_2+1}$, $U_{-l_3^e,l_2^e+1,l_1^e+1,-l_0^e}$), we have $\tilde{A} \phi (x)=0$ for $\phi (x) \in U_{-l_0 ,-l_1,l_2+1 ,l_3+1}$ (resp. $\phi (x) \in  U_{-l_0,l_1+1,-l_2,l_3+1} $, $\phi (x) \in U_{-l_0 ,l_1+1,l_2+1,-l_3}$).
\end{proof}
It follows from Proposition \ref{prop:tAV} that the kernel of the operator $\tilde {A}$ coincides with the space $V$. We denote the monic characteristic polynomial of the operator $H^{(l_0,l_1,l_2, l_3)}$ on the space $V$ by $P(E)$.
For simplicity, we set $u(x)=  \sum_{i=0}^3 l_i(l_i+1)\wp (x+\omega_i)$ and $H=H^{(l_0,l_1,l_2, l_3)}= -d^2/dx^2 +u(x)$.
\begin{prop}
Set $\tilde{a}_0(x)=1$ and $\tilde{a}_{g+1}(x)=0$. The operator $\tilde {A}$ is written as 
\begin{equation}
\tilde {A}= (-1)^g\left[ \sum_{j=0}^{g} \left\{ \tilde{a}_j(x)\frac{d}{dx}-\frac{1}{2} \left( \frac{d}{dx} \tilde{a}_j(x) \right) \right\} H ^{g-j} + \sum_{j=0}^{g} c_j H^{g-j} \right] , \label{Adef0}
\end{equation}
for even doubly-periodic functions $\tilde{a}_j(x)$ $(j=1,\dots ,g)$ and constants $c_j$ $(j=0,\dots ,g)$, and $\tilde{a}_j(x)$ $(j=0,\dots ,g)$ satisfies
\begin{equation}
\tilde{a}'''_j(x)-4u(x)\tilde{a}'_j(x)+4\tilde{a}'_{j+1}(x)-2u'(x)\tilde{a}_j(x)=0. \label{eq:a'''}
\end{equation}
\end{prop}
\begin{proof}
Since $\tilde {A}$ is a monic differential operator of order $2g+1$, it is written as 
\begin{equation*}
\tilde {A}= (-1)^g\left[ \sum_{j=0}^{g} \left( \tilde{a}_j(x)\frac{d}{dx} +\tilde{b}_j(x) \right) H ^{g-j}\right] ,
\end{equation*}
where $\tilde{a}_0(x)=1$.
We have
\begin{align*}
& [(-1)^g \tilde {A},H]= \sum_{j=0}^{g} \left[ \tilde{a}_j(x)\frac{d}{dx} +\tilde{b}_j(x), - \frac{d^2}{dx^2} +u(x) \right] H^{g-j} \\
& =  \sum_{j=0}^{g} \left( \tilde{a}_j(x)u'(x)+2\tilde{a}'_j(x)\frac{d^2}{dx^2}+(\tilde{a}_j''(x)+2\tilde{b}'_j (x)) \frac{d}{dx}+\tilde{b}''_j (x)\right)H^{g-j} \nonumber \\
& = \sum_{j=0}^{g} \left(2\tilde{a}'_j(x)(-H+u(x))+(\tilde{a}_j''(x)+2\tilde{b}'_j (x)) \frac{d}{dx}+ \tilde{a}_j(x)u'(x)+\tilde{b}''_j (x)\right)H^{g-j} \nonumber \\
& = \sum_{j=0}^{g} \left((\tilde{a}_j''(x)+2\tilde{b}'_j (x)) \frac{d}{dx}-2\tilde{a}'_{j+1} (x)+ 2\tilde{a}'_j(x)u(x)+ \tilde{a}_j(x)u'(x)+\tilde{b}''_j (x) \right)H^{g-j} \nonumber \\
& =0. \nonumber
\end{align*}
Hence we have 
\begin{align*}
& \tilde{a}_j''(x)+2\tilde{b}'_j (x)=0, \quad -2\tilde{a}'_{j+1} (x)+ 2\tilde{a}'_j(x)u(x)+ \tilde{a}_j(x)u'(x)+\tilde{b}''_j (x) =0.
\end{align*}
Therefore
\begin{align*}
& \tilde{b}_j(x)= -\tilde{a}'_j (x)/2 +c_j, \quad \tilde{a}'''_j(x)-4u(x)\tilde{a}'_j(x)+4\tilde{a}'_{j+1}(x)-2u'(x)\tilde{a}_j(x)=0
\end{align*}
for some constants $c_j$ $(j=0,\dots ,g)$, and we obtain the proposition.
\end{proof}

Set
\begin{equation}
\tilde{ \Xi}(x,E) = \sum_{i=0}^{g} \tilde{a}_{g-i}(x) E^i. \label{Xiag}
\end{equation}
It follows from Eq.(\ref{eq:a'''}) that
\begin{align}
& \left( \frac{d^3}{dx^3}-4\left( u(x)-E\right)\frac{d}{dx}-2u'(x)  \right) \tilde{ \Xi}(x,E)=0.
\label{prodDE} 
\end{align}
Set
\begin{align}
 & \tilde{Q}(E)=  \tilde{ \Xi} (x,E)^2\left( E- u(x)\right) +\frac{1}{2}\tilde{ \Xi} (x,E)\frac{d^2\tilde{ \Xi} (x,E)}{dx^2}-\frac{1}{4}\left(\frac{d\tilde{ \Xi} (x,E)}{dx} \right)^2. \label{const}
\end{align}
It is shown by differentiating the right hand side of Eq.(\ref{const}) and applying Eq.(\ref{prodDE}) that $\tilde{Q}(E)$ is independent of $x$. $\tilde{Q}(E)$ is a monic polynomial in $E$ of degree $2g+1$, which follows from the expression for $\tilde{\Xi }(x,E)$ given by Eq.(\ref{Xiag}).
Similarly to Proposition 3.2 in \cite{Tak3}, we can show 
\begin{equation*}
\left((-1)^g \tilde {A} - \sum_{j=0}^{g} c_j H^{g-j}\right)^2=\left(\sum_{j=0}^{g} \left\{ \tilde{a}_j(x)\frac{d}{dx}-\frac{1}{2} \left( \frac{d}{dx} \tilde{a}_j(x) \right) \right\} H ^{g-j}\right)^2= \tilde{Q}(H)^2. 
\end{equation*}

Let us recall the function $\Xi (x,E)$ defined in \cite{Tak1}.
It satisfies Eq.(\ref{prodDE}) and has an expression 
\begin{equation}
\Xi (x,E)=c_0(E)+\sum_{i=0}^3 \sum_{j=0}^{l_i-1} b^{(i)}_j (E)\wp (x+\omega_i)^{l_i-j},
\label{Fx}
\end{equation}
where the coefficients $c_0(E)$ and $b^{(i)}_j(E)$ are polynomials in $E$, they do not have common divisors and the polynomial $c_0(E)$ is monic.
It is shown that the dimension of the functions which are doubly-periodic and satisfy Eq.(\ref{prodDE}) is one (see \cite[Proposition 3.9]{TakP}).
Since $\tilde{\Xi }(x,E)$ is a polynomial with respect to the variable $E$ and coefficients of $\Xi (x,E)$ are coprime, we have 
\begin{equation}
\tilde{\Xi }(x,E) = \Xi (x,E) \left( \sum _{i=0}^k d _i E^{k-i}\right) \label{eq:tXiXi}
\end{equation}
 for non-negative integer $k$ and constants $d_i$ ($i=0,\dots ,k$) such that $d_0 =1$.
Write 
\begin{equation}
\Xi (x,E) = \sum_{i=0}^{g-k} a_{g-k-i}(x) E^i. \label{Xiag0}
\end{equation}
Then we have $a_0(x)=1$ and
\begin{equation*}
\tilde{a}_{g-i} (x) = \sum_{j=0}^{k} a_{g-i-k+j}(x) d_{k-j}.
\end{equation*}
The functions $a_i(x)$ ($i=0, \dots ,g-k$) also satisfy Eq.(\ref{eq:a'''}), because $\Xi (x,E) $ satisfies Eq.(\ref{prodDE}).
Set 
\begin{equation}
A= \sum_{j=0}^{g-k} \left\{ a_j(x)\frac{d}{dx}-\frac{1}{2} \left( \frac{d}{dx} a_j(x) \right) \right\} H ^{g-k-j}, \label{Adef00}
\end{equation}
It follows from Eq.(\ref{eq:a'''}) for $a_i (x)$ that $[A,H]=0$.
\begin{prop} \label{prop:tAA}
\begin{equation}
(-1)^g \tilde{A}- \sum_{l=0}^{g} c_{g-l} H^{l}=A \left( \sum _{j=0}^k d _{k-j} H^{j}\right). \label{eq:AAtil}
\end{equation}
\end{prop}
\begin{proof}
\begin{align*}
& (-1)^g \tilde{A}- \sum_{l=0}^{g} c_{g-l} H^{l}= \sum_{i=0}^{g} \left\{ \tilde{a}_i(x)\frac{d}{dx}-\frac{1}{2} \left( \frac{d}{dx} \tilde{a}_i(x) \right) \right\} H^{g-i} \\
& = \sum_{i=0}^{g} \left\{  \left(\sum_{j=0}^{k} a_{g-i-k+j}(x) d_{k-j}\right)\frac{d}{dx}-\frac{1}{2} \left( \frac{d}{dx} \sum_{j=0}^{k} a_{g-i-k+j}(x) d_{k-j} \right) \right\} H^{g-i} \nonumber \\
& = \sum_{i=0}^{g} \sum_{j=0}^{k} \left\{  a_{g-i-k+j}(x) \frac{d}{dx}-\frac{1}{2} \left( \frac{d}{dx}  a_{g-i-k+j}(x) \right) \right\} H^{g-i-j} d_{k-j} H^j \nonumber \\
& =A \left(\sum _{j=0}^k d _{k-j} H^j\right) .\nonumber
\end{align*}
\end{proof}
Set
\begin{align}
 & Q(E)= \Xi (x,E)^2\left( E- u(x)\right) +\frac{1}{2} \Xi (x,E)\frac{d^2 \Xi (x,E)}{dx^2}-\frac{1}{4}\left(\frac{d \Xi (x,E)}{dx} \right)^2. \label{const0}
\end{align}
Then the right hand side of Eq.(\ref{const0}) is independent of $x$, and $Q(E)$ is a monic polynomial in $E$ of degree $2(g-k)+1$.
It is shown in \cite[Proposition 3.2]{Tak3} that $A^2=Q(H)$.
By Eq.(\ref{eq:tXiXi}) and definitions of $Q(E)$ and $\tilde{Q}(E) $, we have 
\begin{equation}
\tilde{Q}(E) =Q(E) \left( \sum _{i=0}^k d _i E^{k-i}\right)^2, \quad 
\deg Q(E) \leq \deg \tilde{Q} (E) =2g+1.
\label{eq:degQ}
\end{equation}
On zeros of $Q(E)$ and $P(E)$, the following proposition is shown in \cite[Theorem 3.8]{Tak1} (see also \cite[Proposition 2.4]{Tak3}):
\begin{prop} (c.f. \cite[Theorem 3.8]{Tak1}) \label{prop:zeros}
The set of zeros of $Q(E)$ coincides with the set of zeros of $P(E)$.
\end{prop}
\begin{prop} \label{thm:dist} (c.f. \cite{TakR})
Roots of the equation $P(E)=0$ are distinct for generic periods $(2\omega_1, 2\omega _3)$.
\end{prop}
We prove Proposition \ref{thm:dist} in the appendix. 
\begin{prop}
For the periods $(2\omega_1, 2\omega _3)$ that roots of the equation $P(E)=0$ are distinct, we have $P(E)=Q(E)=\tilde{Q}(E)$ and $\Xi (x,E)=\tilde{\Xi }(x,E)$.
\end{prop}
\begin{proof}
By assumption, the equation $P(E)=0$ has $2g+1$ roots which are distinct.
It follows from Proposition \ref{prop:zeros} that the number of roots of the equation $Q(E)=0$ is equal to or more than $2g+1$. Then we have $\deg Q(E) \geq 2g+1$. Combining with Eq.(\ref{eq:degQ}) and $\deg P(E)=2g+1$, we have $\tilde{Q}(E) =Q(E)=P(E)$, and the value $k$ in Eq.(\ref{eq:tXiXi}) is zero. Hence $\Xi (x,E)=\tilde{\Xi }(x,E)$.
\end{proof}
\begin{prop} \label{prop:cj0}
All constants $c_j$ $(j=1,\dots ,g)$ are zero. Namely we have
\begin{equation}
\tilde {A}= (-1)^g \sum_{j=0}^{g} \left\{ \tilde{a}_j(x)\frac{d}{dx}-\frac{1}{2} \left( \frac{d}{dx} \tilde{a}_j(x) \right) \right\} H ^{g-j} . \label{Adef1}
\end{equation}
\end{prop}
\begin{proof}
First, we assume that roots of the equation $P(E)=0$ are distinct. Let $\{E _i\}_{i=1,\dots,2g+1}$ be the roots of the equation $P(E)=0$. Then $E _i$ $(i=1,\dots,2g+1)$ is an eigenvalues of the operator $H$ on the space $V$. Let $\phi _i(x) \in V$ be the eigenfunction of the operator $H$ with the eigenvalue $E_i$. It follows from Proposition \ref{prop:tAA} that $((-1)^g \tilde{A}- \sum_{l=0}^{g} c_{g-l} H^{l})^2=A^2=Q(H)=P(H)$. Hence $ (\tilde{A} -2(-1)^g \sum_{l=0}^{g} c_{g-l} H^{l}) \tilde{A} +(\sum_{l=0}^{g} c_{g-l} H^{l})^2=P(H)$. We apply this operator to $\phi _i(x)$. Since the operator $\tilde{A}$ annihilates the space $V$, we have
$(\sum_{l=0}^{g} c_{g-l} (E_i)^{l})^2 \phi _i(x)= P(E_i) \phi _i(x) =0$. Hence $\sum_{l=0}^{g} c_{g-l} (E_i)^{l} =0$ for $i=1,\dots ,g$. Since the polynomial of degree less that $g+1$ cannot have $2g+1$ zeros, we have $c_i=0$ $(i=0,\dots ,g)$ for the case that roots of the equation $P(E)=0$ are distinct.

Now consider the case that roots of the equation $P(E)=0$ are not distinct. Since the operator $\tilde{A}$ is defined by Eq.(\ref{eq:Ate}) or Eq.(\ref{eq:Ato}), the coefficients $c_i$ $(i=0,\dots ,g)$ are continuous as a function of periods $(2\omega_1, 2\omega _3)$. Since $c_i=0$ $(i=0,\dots ,g)$ for a dense set of periods, we have $c_i=0$ $(i=0,\dots ,g)$ for all periods.
\end{proof}

For the case that roots of the equation $P(E)=0$ are distinct, it is shown that $P(E)=Q(E)$, the value $k$ in Eq.(\ref{eq:tXiXi}) is zero, $\Xi (x,E)=\tilde{\Xi }(x,E)$ and $\tilde{A}=(-1)^g A$.

We now consider the case that roots of the equation $P(E)=0$ are not distinct.
It is already shown that $P(E)= \tilde{Q}(E)$ for a dense set of periods, $\deg P(E)=\deg \tilde{Q}(E)$ for all periods, and coefficients of $P(E)$ and $Q(E)$ are continuous with respect to the periods. Hence we have $P(E)= \tilde{Q}(E)$ for all periods. 
Assume that the value $k$ in Eq.(\ref{eq:tXiXi}) is positive.
Combining with Eq.(\ref{eq:degQ}) and Proposition \ref{prop:zeros}, all zeros of $(\sum _{i=0}^k d _i E^i)^2$ are zeros of $Q(E)$. Let $E_0$ be a zero of $\sum _{i=0}^k d _i E^i$, i.e. $\sum _{i=0}^k d _i (E_0)^i=0$.
By Eq.(\ref{eq:AAtil}) and Proposition \ref{prop:cj0}, it is shown that, if $f(x)$ satisfies $Hf(x)=E_0f(x)$, then $\tilde{A} f(x)=0$. Hence all solutions to $(H-E_0)f(x) =0$ are contained in the space $V$. But it contradicts that all solutions to $(H-E_0)f(x) =0$ cannot be doubly-periodic up to signs, which is essentially shown in \cite[Theorem 3.8]{Tak1} (see also the proof of \cite[Proposition 3.9]{TakP}).
Hence we have $k=0$. Therefore $\tilde{A} =(-1)^g A$, $\deg P(E)= \deg Q(E)$ and $P(E)=Q(E)$ for all periods.
It follows from Proposition \ref{prop:tAV} that the operator $(-1)^g A(=\tilde{A})$ is characterized by the monic operator of order $2g+1$ which annihilates all elements in the space $V$. 
Summarizing, we obtain the following theorem:
\begin{thm} \label{thm:fing}
Let $\tilde{A}$ be the operator defined by composing generalized Darboux transformation (see Eq.(\ref{eq:Ate}) or Eq.(\ref{eq:Ato})) and $A$ be the operator defined from the even doubly-periodic function $\Xi (x,E)$ (see Eqs.(\ref{Xiag0}, \ref{Adef00})). Let $P(E)$ be the characteristic polynomial of the operator $H$ on the space $V$, and $Q(E)$ be the polynomial defined from $\Xi (x,E)$ (see Eq.(\ref{const0})).
\\
(i) We have $\tilde{A} =(-1)^g A$ and $P(E)=Q(E)$.\\
(ii) The operator $(-1)^g A$ is also characterized by the monic operator of order $2g+1$ which annihilates all elements in the space $V$.
\end{thm}
Note that we proved Conjecture 1 in \cite{Tak3} by (i) and Conjecture 2 in \cite{Tak3} by (ii).
Since $P(E)=Q(E)$, the genus of the spectral curve $\nu^2 =-Q(E)$ is $g$, where $g$ is defined in Eq.(\ref{eq:ge}) ($l_0+l_1+l_2+l_3$: even) or Eq.(\ref{eq:go}) ($l_0+l_1+l_2+l_3$: odd), and it agrees with the results in \cite{GW,Tre}.

We consider isomonodromic property again. Let $l_i ^e$, $l_i ^o$ $(i=0,1,2,3)$ be the numbers defined in Eqs.(\ref{def:lie}, \ref{def:lio}). We have shown in section \ref{sec:int} that, if $l_0+l_1+l_2+l_3$ is even (resp. odd), then the operator $H^{(l_0,l_1,l_2,l_3)}$ is linked to $H^{(\tilde{l}_0,\tilde{l}_1,\tilde{l}_2,\tilde{l}_3)}$ by generalized Darboux transformation, where $H^{(\tilde{l}_0,\tilde{l}_1,\tilde{l}_2,\tilde{l}_3)}$ is any operator listed in Eq.(\ref{eq:gDHe}) (resp. Eq.(\ref{eq:gDHo})). We now show that functions related to monodromy for the operator $H^{(l_0,l_1,l_2,l_3)}$ coincides to ones for the operator $H^{(\tilde{l}_0,\tilde{l}_1,\tilde{l}_2,\tilde{l}_3)}$. For this purpose, we recall functions defined in \cite{Tak1,Tak3,Tak4}.

Let $\Xi ^{(l_0,l_1,l_2,l_3)}(x,E)$ be the function defined in Eq.(\ref{Fx}), $Q ^{(l_0,l_1,l_2,l_3)}(E)$ be the polynomial defined in Eq.(\ref{const0}) and $P ^{(l_0,l_1,l_2,l_3)}(E)$ be the polynomial $P(E)$ defined in this section.
Set 
\begin{equation*}
\Lambda ^{(l_0,l_1,l_2,l_3)} (x,E)=\sqrt{\Xi ^{(l_0,l_1,l_2,l_3)}(x,E)}\exp \int \frac{\sqrt{-Q^{(l_0,l_1,l_2,l_3)}(E)}dx}{\Xi ^{(l_0,l_1,l_2,l_3)}(x,E)}.
\end{equation*}
Then the function $\Lambda ^{(l_0,l_1,l_2,l_3)} (x,E)$ is a solution to the differential equation
\begin{equation*}
(H^{(l_0,l_1,l_2,l_3)}-E)f(x)=0.
\end{equation*}
Set
\begin{equation*}
\Phi _i(x,\alpha )= \frac{\sigma (x+\omega _i -\alpha ) }{ \sigma (x+\omega _i )} \exp (\zeta( \alpha )x), \quad \quad (i=0,1,2,3),
\end{equation*}
where $\sigma (x)$ is the Weierstrass sigma function.
The function $\Lambda ^{(l_0,l_1,l_2,l_3)} (x,E)$ admits an expression in terms of Hermite-Krichever Ansatz. Namely, we have the following proposition:
\begin{prop} (c.f. \cite[Theorem 2.3]{Tak4}) \label{thm:alpha}
The function $\Lambda ^{(l_0,l_1,l_2,l_3)} (x,E)$ is expressed as 
\begin{align}
& \Lambda ^{(l_0,l_1,l_2,l_3)} (x,E) = \exp \left( \kappa x \right) \left( \sum _{i=0}^3 \sum_{j=0}^{l_i-1} \tilde{b} ^{(i)}_j \left( \frac{d}{dx} \right) ^{j} \Phi _i(x, \alpha ) \right)
\label{Lalpha}
\end{align}
for some values $\tilde{b} ^{(i)}_j$ $(i=0,\dots ,3, \: j= 0,\dots ,l_i-1)$ except for finitely-many eigenvalues $E$. The values $\alpha $ and $\kappa $ are expressed as
\begin{align*}
&  \wp (\alpha ) =R^{(l_0,l_1,l_2,l_3)} _1(E), \; \; \; \wp ' (\alpha ) = R^{(l_0,l_1,l_2,l_3)} _2 (E)\sqrt{-Q^{(l_0,l_1,l_2,l_3)}(E)} , \\
& \; \; \kappa  =R_3^{(l_0,l_1,l_2,l_3)}  (E) \sqrt{-Q^{(l_0,l_1,l_2,l_3)}(E)}, \nonumber
\end{align*}
where $R_1 ^{(l_0,l_1,l_2,l_3)} (E) ,R_2 ^{(l_0,l_1,l_2,l_3)} (E) ,R_3^{(l_0,l_1,l_2,l_3)} (E)$ are rational functions in $E$.
\end{prop}
It follows from Eq.(\ref{Lalpha}) that
\begin{align} 
& \Lambda ^{(l_0,l_1,l_2,l_3)} (x+2\omega _k,E)  = \exp (-2\eta _k \alpha +2\omega _k \zeta (\alpha ) +2 \kappa \omega _k ) \Lambda ^{(l_0,l_1,l_2,l_3)} (x,E)   \label{ellint}
\end{align}
for $k=1,2,3$. 
A monodromy formula in terms of hyperelliptic integral was obtained in \cite{Tak3}.
\begin{prop} (c.f. \cite[Theorem 3.7]{Tak3}) \label{thm:conj3} 
 Let  $k \in \{ 1,2,3\}$, $q_k \in \{0,1\}$ and $E_0$ be the eigenvalue such that $\Lambda ^{(l_0,l_1,l_2,l_3)}(x+2\omega _k,E_0)=(-1)^{q_k} \Lambda ^{(l_0,l_1,l_2,l_3)}(x,E_0)$. Then
\begin{align}
& \Lambda ^{(l_0,l_1,l_2,l_3)}(x+2\omega _k,E)=(-1)^{q_k} \Lambda ^{(l_0,l_1,l_2,l_3)}(x,E) \cdot  \label{hypellint} \\
& \quad \quad \quad \exp \left( -\frac{1}{2} \int_{E_0}^{E}\frac{ -2\eta _k a^{(l_0,l_1,l_2,l_3)}(\tilde{E}) +2\omega _k c^{(l_0,l_1,l_2,l_3)}(\tilde{E}) }{\sqrt{-Q^{(l_0,l_1,l_2,l_3)}(\tilde{E})}} d\tilde{E}\right) , \nonumber
\end{align}
where $a^{(l_0,l_1,l_2,l_3)}(E)$ (resp. $c^{(l_0,l_1,l_2,l_3)}(E)$) is a polynomial defined in \cite{Tak4}.
\end{prop}
The following proposition states the coincidence of functions.
\begin{prop} \label{prop:PQRac}
Let $H^{(\tilde{l}_0,\tilde{l}_1,\tilde{l}_2,\tilde{l}_3)}$ be any operator listed in Eq.(\ref{eq:gDHe}) ($l_0+l_1+l_2+l_3$: even) or Eq.(\ref{eq:gDHo}) ($l_0+l_1+l_2+l_3$: odd). Then we have 
\begin{align*}
& P^{(l_0,l_1,l_2,l_3)} (E) = P^{(\tilde{l}_0,\tilde{l}_1,\tilde{l}_2,\tilde{l}_3)}(E) , \quad Q^{(l_0,l_1,l_2,l_3)}(E)  = Q^{(\tilde{l}_0,\tilde{l}_1,\tilde{l}_2,\tilde{l}_3)}(E) , \\
& R_i^{(l_0,l_1,l_2,l_3)}(E)  = R_i^{(\tilde{l}_0,\tilde{l}_1,\tilde{l}_2,\tilde{l}_3)}(E) 
,\quad \quad (i=1,2,3), \nonumber \\
& a^{(l_0,l_1,l_2,l_3)}(E)  = a^{(\tilde{l}_0,\tilde{l}_1,\tilde{l}_2,\tilde{l}_3)}(E) , \quad c^{(l_0,l_1,l_2,l_3)}(E)  = c^{(\tilde{l}_0,\tilde{l}_1,\tilde{l}_2,\tilde{l}_3)}(E) . \nonumber
\end{align*}
\end{prop}
\begin{proof}
We show the proposition for the case that $l_0+l_1+l_2+l_3$ is even
 and $(\tilde{l}_0,\tilde{l}_1,\tilde{l}_2,\tilde{l}_3) =(l^e_3,l^e_2,l^e_1,l^e_0)$. For the other cases, it is shown similarly.

The space $V$ for the operator $H^{(l_0,l_1,l_2,l_3)}$ is written as
\begin{equation*}
U_{-l_0,-l_1,-l_2,-l_3}\oplus U_{-l_0 ,-l_1,l_2+1 ,l_3+1}\oplus  U_{-l_0,l_1+1,-l_2,l_3+1}\oplus  U_{-l_0 ,l_1+1,l_2+1,-l_3},
\end{equation*}
and the corresponding space for the operator $H^{(l^e_3,l^e_2,l^e_1,l^e_0)}$ is written as
\begin{equation*}
U_{-l^e_3,-l^e_2,-l^e_1,-l^e_0}\oplus U_{-l^e_3 ,-l^e_2,l^e_1+1 ,l^e_0+1}\oplus  U_{-l^e_3,l^e_2+1,-l^e_1,l^e_0+1}\oplus  U_{l^e_3 +1,-l^e_2,-l^e_1,l^e_0 +1}.
\end{equation*}
It is seen that the space $U_{-l_0,-l_1,-l_2,-l_3}$ (resp. $U_{-l_0 ,-l_1,l_2+1 ,l_3+1}$, $U_{-l_0,l_1+1,-l_2,l_3+1}$, $U_{-l_0 ,l_1+1,l_2+1,-l_3}$) is linked to the space $U_{-l^e_3,-l^e_2,-l^e_1,-l^e_0}$ (resp. $U_{-l^e_3 ,-l^e_2,l^e_1+1 ,l^e_0+1}$, $U_{-l^e_3,l^e_2+1,-l^e_1,l^e_0+1}$, $U_{l^e_3 +1,-l^e_2,-l^e_1,l^e_0 +1}$) by the generalized Darboux transformation $\tilde{L}_{-l_0,-l_1,-l_2,-l_3}$ (resp. $\tilde{L}_{-l_0 ,-l_1,l_2+1 ,l_3+1}$, $\tilde{L}_{-l_0,l_1+1,-l_2,l_3+1}$, $\tilde{L}_{-l_0 ,l_1+1,l_2+1,-l_3}$) and the shift $x \rightarrow x+\omega _3$ (resp. $x \rightarrow x+\omega _1$, $x \rightarrow x+\omega _2$, $x \rightarrow x$). It follows from Proposition \ref{prop:HL0dL0dH} that the characteristic polynomial of $H^{(l_0,l_1,l_2,l_3)}$ on the space $U_{-l_0,-l_1,-l_2,-l_3}$ (resp. $U_{-l_0 ,-l_1,l_2+1 ,l_3+1}$, $U_{-l_0,l_1+1,-l_2,l_3+1}$, $U_{-l_0 ,l_1+1,l_2+1,-l_3}$) is equal to that of $H^{(l^e_3,l^e_2,l^e_1,l^e_0)}$ on the space $U_{-l^e_3,-l^e_2,-l^e_1,-l^e_0}$ (resp. $U_{-l^e_3 ,-l^e_2,l^e_1+1 ,l^e_0+1}$, $U_{-l^e_3,l^e_2+1,-l^e_1,l^e_0+1}$, $U_{l^e_3 +1,-l^e_2,-l^e_1,l^e_0 +1}$).
Since the polynomial $P(E)$ is written as the product of characteristic polynomials of invariant subspaces, we have $P^{(l_0,l_1,l_2,l_3)}(E)=P^{(l^e_3,l^e_2,l^e_1,l^e_0)}(E)$. It follows from Theorem \ref{thm:fing} that $Q^{(l_0,l_1,l_2,l_3)}(E)=Q^{(l^e_3,l^e_2,l^e_1,l^e_0)}(E)$.

Let $\tilde{\alpha }$ and $\tilde{\kappa }$ be the corresponding values in Eq.(\ref{ellint}) for the parameters $(l^e_3,l^e_2,l^e_1,l^e_0)$. Namely, they satisfy
\begin{align*} 
& \Lambda ^{(l^e_3,l^e_2,l^e_1,l^e_0)} (x+2\omega _k,E)  = \exp (-2\eta _k \tilde{\alpha } +2\omega _k \zeta (\tilde{\alpha } ) +2 \tilde{\kappa } \omega _k ) \Lambda ^{(l^e_3,l^e_2,l^e_1,l^e_0)} (x,E) 
\end{align*}
for $k=1,2,3$.
The monodromy matrix is preserved by the generalized Darboux transformation for almost eigenvalues $E$. Since the values $\exp (\pm (-2\eta _k \alpha +2\omega _k \zeta (\alpha ) +2 \kappa \omega _k ))$ (resp. $\exp (\pm (-2\eta _k \tilde{\alpha } +2\omega _k \zeta (\tilde{\alpha } ) +2 \tilde{\kappa } \omega _k )) $) are eigenvalues of the monodromy matrix, we have
\begin{align*}
& -2\eta _1 \alpha  +2\omega _1 (\zeta (\alpha ) + \kappa ) =\pm (-2\eta _1 \tilde{\alpha }  +2\omega _1 (\zeta (\tilde{\alpha }) + \tilde{\kappa }))+ 2n_1 \pi \sqrt{-1}, \\
& -2\eta _3 \alpha  +2\omega _3 (\zeta (\alpha ) + \kappa ) =\pm (-2\eta _3 \tilde{\alpha }  +2\omega _3 (\zeta (\tilde{\alpha }) + \tilde{\kappa }))+ 2n_3 \pi \sqrt{-1},
\end{align*}
for some integers $n_1$ and $n_3$ and almost eigenvalues $E$.
By the Legendre's relation $\eta _1 \omega _3 - \eta _3 \omega _1 =\pi\sqrt{-1}/2$ and the relation $\zeta (x+2\omega _k)= \zeta(x) +2\eta _k$ $(k=1,3)$, it follows that
\begin{align*}
& \alpha = \pm \tilde{\alpha } -( 2n_1 \omega _3 - 2n_3 \omega _1), \quad \kappa =\pm \tilde{\kappa }.
\end{align*}
It follows from the asymptotic of $\kappa $ (see \cite[Proposition 3.2]{Tak4}) that the sign $\pm$ is plus.
Hence $\wp (\alpha )= \wp (\tilde{\alpha })$, $\wp '(\alpha )=  \wp '(\tilde{\alpha })$ and $\kappa = \tilde{\kappa }$ for almost $E$.
Since $R_i^{(l_0,l_1,l_2,l_3)}(E)$ and $R_i^{(l^e_3,l^e_2,l^e_1,l^e_0)}(E) $ $(i=1,2,3)$ are rational functions, we have $R_i^{(l_0,l_1,l_2,l_3)}(E)  = R_i^{(l^e_3,l^e_2,l^e_1,l^e_0)}(E) $ for $i=1,2,3$.

Let $\tilde{E}_0$ be the corresponding values in Eq.(\ref{hypellint}) for the parameters $(l^e_3,l^e_2,l^e_1,l^e_0)$.
By applying a similar discussion for Eq.(\ref{hypellint}), we obtain that the integrals
\begin{align*}
& \int_{E_0}^{E}\frac{ a^{(l_0,l_1,l_2,l_3)}(\tilde{E})}{\sqrt{-Q^{(l_0,l_1,l_2,l_3)}(\tilde{E})}} d\tilde{E} -\int_{\tilde{E}_0}^{E}\frac{ a^{(l^e_3,l^e_2,l^e_1,l^e_0)}(\tilde{E})}{\sqrt{-Q^{(l^e_3,l^e_2,l^e_1,l^e_0)}(\tilde{E})}} d\tilde{E}, \\
& \int_{E_0}^{E}\frac{ c^{(l_0,l_1,l_2,l_3)}(\tilde{E})}{\sqrt{-Q^{(l_0,l_1,l_2,l_3)}(\tilde{E})}} d\tilde{E} -\int_{\tilde{E}_0}^{E}\frac{ c^{(l^e_3,l^e_2,l^e_1,l^e_0)}(\tilde{E})}{\sqrt{-Q^{(l^e_3,l^e_2,l^e_1,l^e_0)}(\tilde{E})}} d\tilde{E}, \\
\end{align*}
are constants. By differentiating them in the variable $E$ and using the relation $Q^{(l_0,l_1,l_2,l_3)}(E)  = Q^{(l^e_3,l^e_2,l^e_1,l^e_0)}(E)$, we have $a^{(l_0,l_1,l_2,l_3)}(E)  = a^{(l^e_3,l^e_2,l^e_1,l^e_0)}(E)$ and $c^{(l_0,l_1,l_2,l_3)}(E)  = c^{(l^e_3,l^e_2,l^e_1,l^e_0)}(E)$.
\end{proof}
It follows from Proposition \ref{prop:PQRac} that the hyperelliptic-to-elliptic integral formulae obtained in \cite{Tak4} for the parameters $(l_0,l_1,l_2,l_3)$ coincide with that for the parameters $(\tilde{l}_0,\tilde{l}_1,\tilde{l}_2,\tilde{l}_3)$.

We calculate the commuting operator $A$ explicitly for the case $l_0=l_1=l_2=l_3=g \in \Zint _{\geq 1}$.
The commuting operator $A$ of $H^{(g,g,g,g)}$ is written as
\begin{equation*}
A= (-1)^g L_{-g,-g,-g,-g}= (-1)^g \wp '(x)^{2g+1} \widehat{\Phi}(\wp (x)) \circ \left(\frac{1}{\wp' (x)}\frac{d}{dx} \right)^{2g+1} \circ \widehat{\Phi}(\wp (x)) ^{-1},
\end{equation*}
(see Eq.(\ref{op:L})), and $\widehat{\Phi}(\wp (x))=((\wp(x) -e_1)^{-g/2}(\wp(x)-e_2)^{-g/2}(\wp(x)-e_3))^{-g/2}= (\wp '(x)/2)^{-g}$. Hence we have
\begin{equation}
A= (-1)^g \wp '(x)^{g+1} \circ \left(\frac{1}{\wp' (x)}\frac{d}{dx} \right)^{2g+1} \circ \wp '(x)^{g}.
\label{op:AHgggg}
\end{equation}
On the other hand, we have $\sum _{i=0}^3 \wp (x+\omega _i) = 4\wp (2x)$. By the change $2x \rightarrow x$, we recover the Lam\'e operator $4(-d^2/dx^2 +g(g+1)\wp (x))$. Hence the commuting operator $A$ of $H^{(g,0,0,0)}$ is written as 
\begin{equation}
A= (-1)^g \wp '(x/2)^{g+1} \circ \left(\frac{1}{\wp' (x/2)}\frac{d}{dx} \right)^{2g+1} \circ \wp '(x/2)^{g}.
\label{op:AHg000}
\end{equation}
Therefore we obtain the following proposition:
\begin{prop}
The commuting operator $A$ for the Lam\'e operator $H^{(g,0,0,0)}$ $(g \in \Zint _{\geq 1})$ is written as Eq.(\ref{op:AHg000}).
\end{prop}

\appendix
\section {Elliptic functions} \label{sect:append}
This appendix presents the definitions of and the formulas for the elliptic functions.
The Weierstrass $\wp$-function is defined by
\begin{align}
& \wp (x)= \frac{1}{x^2}+  \sum_{(m,n)\in \Zint \times \Zint \setminus \{ (0,0)\} } \left( \frac{1}{(x-2m\omega_1 -2n\omega_3)^2}-\frac{1}{(2m\omega_1 +2n\omega_3)^2}\right).
\end{align}
Setting $\omega_2=-\omega_1-\omega_3$ and $e_k=\wp(\omega_k) $ $(k=1,2,3)$
yields the relations
\begin{align}
& e_1+e_2+e_3=0, \; \; \; \label{eq:Leg} \\
& \frac{\wp''(x)}{(\wp'(x))^2}=\frac{1}{2}\left( \frac{1}{\wp (x)-e_1}+\frac{1}{\wp (x)-e_2}+\frac{1}{\wp(x)-e_3} \right), \nonumber \\
& \wp(x+\omega_i)=e_i+\frac{(e_i-e_{i'})(e_i-e_{i''})}{\wp(x)-e_i} \; \; \; \; (i=1,2,3),\nonumber
\end{align}
where $i', i'' \in \{1,2,3\}$ with $i'<i''$, $i\neq i'$ and $i\neq i''$.

\section{Proofs of Propositions \ref{prop:HL0dL0dH} (ii), \ref{prop:distrt} and \ref{thm:dist}}

Let $\alpha _i$ be a number such that $\alpha _i= -l_i$ or $\alpha _i= l_i+1$ for each $i\in \{ 0,1,2,3\} $.
Set ${\bf a} = (\alpha _0, \alpha _1, \alpha _2, \alpha _3)$ and
\begin{align*}
& v_r^{{\bf a}}= (\wp(x)-e_1)^{\alpha _1/2}(\wp(x)-e_3)^{\alpha _3/2}(\wp(x)-e_2)^{\alpha _2/2+r}, \\
& a_{r+1,r}^{{\bf a}}=-4(r+\gamma_1^{{\bf a}})(r+\gamma_2^{{\bf a}}) , \nonumber \\
& a_{r-1,r}^{{\bf a}}=-4r(r+\alpha_2-1/2)(e_2-e_3)(e_2-e_1), \nonumber \\
& a_{r,r}^{{\bf a}}= -4r((e_2-e_3)(r+\alpha_2+\alpha_1)+(e_2-e_1)(r+\alpha_2+\alpha_3)) \nonumber \\
& \quad \quad \quad \quad -4e_2\gamma_1^{{\bf a}}\gamma_2^{{\bf a}}+e_1(\alpha_2+\alpha_3)^2+e_2(\alpha_1+\alpha_3)^2+e_3(\alpha_1+\alpha_2)^2, \nonumber \\
& \gamma_1^{{\bf a}}=(\alpha _0+ \alpha_1+\alpha_2+\alpha_3)/2, \quad \gamma_2^{{\bf a}}=(-\alpha _0+ \alpha_1+\alpha_2+\alpha_3+1)/2. \nonumber
\end{align*}
Then the action of the operator $H^{(l_0,l_1,l_2,l_3)}$ is written as
\begin{align}
&  H^{(l_0,l_1,l_2,l_3)}  v_r^{{\bf a}} = a_{r+1,r}^{{\bf a}} v_{r+1}^{{\bf a}} +a_{r,r}^{{\bf a}} v_{r} ^{{\bf a}}+a_{r-1,r}^{{\bf a}} v_{r-1}^{{\bf a}} , \label{ze2}
\end{align}
(see \cite{Tak2}). Set $d=-\sum_{i=0}^3 \alpha _i /2$ and assume $d\in \Zint_{\geq 0}$. Then the space $V_{\alpha _0, \alpha _1, \alpha _2, \alpha _3}$ is spanned by $v_0^{{\bf a}}, v_1^{{\bf a}}, \dots , v_d^{{\bf a}}$, and it follows from Eq.(\ref{ze2}) that $H^{(l_0,l_1,l_2, l_3)}$ preserves the space $V_{\alpha _0, \alpha _1, \alpha _2, \alpha _3}$.
The operator $H^{(\alpha _0 +d,\alpha _1 +d,\alpha _2 +d,\alpha _3 +d)}$ also preserves the $d+1$-dimensional space $V_{-\alpha _0-d, -\alpha _1-d, -\alpha _2-d, -\alpha _3-d}$.
\begin{prop} \label{aprop:HL0dL0dH} (Proposition \ref{prop:HL0dL0dH} (ii))
Set $d=-\sum_{i=0}^3 \alpha _i /2$ and assume $d\in \Zint_{\geq 0}$. Then the characteristic polynomial of the operator $H^{(l_0,l_1,l_2,l_3)}$ on the space $V_{\alpha _0, \alpha _1, \alpha _2, \alpha _3}$ coincides with that of the operator $H^{(\alpha _0 +d,\alpha _1 +d,\alpha _2 +d,\alpha _3 +d)}$ on the space $V_{-\alpha _0-d, -\alpha _1-d, -\alpha _2-d, -\alpha _3-d}$.
\end{prop}
\begin{proof}
Set ${\bf -a-d}= (-\alpha _0-d, -\alpha _1-d, -\alpha _2-d, -\alpha _3-d)$. The action of the operator $H^{(\alpha _0+d, \alpha _1+d, \alpha _2+d, \alpha _3+d)}$ on the space $V_{-\alpha _0-d, -\alpha _1-d, -\alpha _2-d, -\alpha _3-d}$ is written as
\begin{align*}
&  H^{(\alpha _0+d, \alpha _1+d, \alpha _2+d, \alpha _3+d)}  v_r^{{\bf -a-d}} = a_{r+1,r}^{{\bf -a-d}} v_{r+1}^{{\bf -a-d}} +a_{r,r}^{{\bf -a-d}} v_{r} ^{{\bf -a-d}}+a_{r-1,r}^{{\bf -a-d}} v_{r-1}^{{\bf -a-d}}. 
\end{align*}
Set $u_r^{{\bf a}}= ((-1)^r d!/(r! (d-r)!)) v_{d-r}^{{\bf a}}$. By a direct calculation, Eq.(\ref{ze2}) is rewritten as
\begin{align*}
&  H^{(l_0,l_1,l_2,l_3)}  u_r^{{\bf a}} = a_{r-1,r}^{{\bf -a-d}} u_{r+1}^{{\bf a}} +a_{r,r}^{{\bf -a-d}} u_{r} ^{{\bf a}}+ a_{r+1,r}^{{\bf -a-d}} u_{r-1}^{{\bf a}} . \end{align*}
Hence the martix representation of $ H^{(\alpha _0+d, \alpha _1+d, \alpha _2+d, \alpha _3+d)}$ on the basis $\langle v_{0}^{{\bf -a-d}} , \dots , v_{d}^{{\bf -a-d}} \rangle$ is the transposed one of $ H^{(l_0,l_1,l_2,l_3)}$ on the basis $\langle u_{0}^{{\bf a}} , \dots , u_{d}^{{\bf a}} \rangle$, and we obtain that the characteristic polynomial of the operator $H^{(l_0,l_1,l_2,l_3)}$ on the space $V_{\alpha _0, \alpha _1, \alpha _2, \alpha _3}$ coincides with that of the operator $H^{(\alpha _0 +d,\alpha _1 +d,\alpha _2 +d,\alpha _3 +d)}$ on the space $V_{-\alpha _0-d, -\alpha _1-d, -\alpha _2-d, -\alpha _3-d}$.
\end{proof}

\begin{prop} (Proposition \ref{prop:distrt}) \label{aprop:distrt}
Set $d=-\sum_{i=0}^3 \alpha_i /2$ and assume $d\in \Zint_{\geq 0}$ and $\alpha _i \neq \alpha _j$ for some $i,j \in \{0,1,2,3\}$. Then zeros of the characteristic polynomial of the operator $H^{(l_0,l_1,l_2,l_3)}$ on the space $V _{\alpha _0, \alpha _1, \alpha _2, \alpha _3}$ are distinct for generic periods $(2\omega_1, 2\omega _3)$.
\end{prop}
\begin{proof}
The assumption  $\alpha _i \neq \alpha_j$ for some $i,j \in \{0,1,2,3\}$ is rewritten that $\alpha _0, \alpha _1 ,\alpha _2, \alpha _3$ do not satisfy $\alpha _0= \alpha _1 =\alpha _2= \alpha _3$, and it is easy to show that there exists $i_0 ,i_1 ,i_2 \in \{0,1,2,3 \}$ such that $i_1 \neq  i_2$, $\alpha _{i_0} \neq \alpha _{i_1}$ and $\alpha _{i_0} \neq \alpha _{i_2}$.
By permutation of periods $\omega _1, \omega _2, \omega _3$ and shifts $x \rightarrow x+\omega _i$ $(i=1,2,3)$, we can permutate numbers $\alpha _0, \alpha _1, \alpha _2, \alpha _3$. Hence it is sufficient to show the proposition under the assumption $\alpha _3 \neq \alpha _1$ and $\alpha _3 \neq \alpha _2$.

Set $\omega _1 =1/2$, $\omega _3 = \tau /2$ and $p=\exp (\pi \sqrt{-1} \tau)$. Then $e_1$, $e_2$ and $e_3$ are expressed as power series in $p$ and we have $e_1= \pi^2 (2/3 +O(p^2))$, $e_2= \pi^2 (-1/3+ 8p +O(p^2))$ and $e_3=\pi^2(-1/3-8 p+O(p^2))$. The matrix elements are expressed as
\begin{align*}
& a_{r+1,r}^{{\bf a}}=\tilde{a}_{r+1,r}^{(0)} , \quad a_{r-1,r}^{{\bf a}}=\tilde{a}_{r-1,r}^{(1)} (p+O(p^2)), \quad  a_{r,r}^{{\bf a}}= \tilde{a}_{r,r}^{(0)} +\tilde{a}_{r,r}^{(1)} p+ O(p^2), \\
& \tilde{a}_{r+1,r}^{(0)} =-4(r+\gamma_1^{{\bf a}})(r+\gamma_2^{{\bf a}}), \quad \tilde{a}_{r-1,r}^{(1)}= 64r(r+\alpha_2-1/2), \nonumber \\
& \textstyle \tilde{a}_{r,r}^{(0)}= (2r+ \alpha_2+\alpha_3)^2- \sum_{i=0}^3 l_i (l_i+1)/3, \nonumber \\
& \tilde{a}_{r,r}^{(1)}= -8\{r(12r+ 8\alpha_1+ 12\alpha_2+4\alpha_3)+4 \gamma_1^{{\bf a}}\gamma_2^{{\bf a}}- (\alpha_1+\alpha_3)^2+(\alpha_1+\alpha_2)^2\}. \nonumber 
\end{align*}
If $p = 0$, then the operator $H^{(l_0,l_1,l_2,l_3)}$ acts triangularly and eigenvalues are written as $\tilde{a}_{r,r}^{(0)}$ ($r=0, \dots ,d$). Since the eigenvalues are quadratic in $r$, multiplicity of the eigenvalues is one or two.
Hence it is sufficient to show that the eigenvalue with multiplicity two on the case $p=0$ separates when $p$ varies.
Assume that $\tilde{a}_{r,r}^{(0)}(= \tilde{a}_{r',r'}^{(0)})$ is the eigenvalue with multiplicity two, $0\leq r<r'\leq d$ and $r,r' \in \Zint$.
Then we have $r+r'=-(\alpha_2+\alpha_3)$.

If $r+1 <r'$, then the eigenvalue around $\tilde{a}_{r,r}^{(0)}$ is expanded as
$E=\tilde{a}_{r,r}^{(0)} + c_1 p + \dots $, and $c_1$ satisfies
\begin{align}
& \{ (\tilde{a}_{r-1,r-1}^{(0)} -\tilde{a}_{r,r}^{(0)})(\tilde{a}_{r,r}^{(1)} -c_1) (\tilde{a}_{r+1,r+1}^{(0)}- \tilde{a}_{r,r}^{(0)}) -(\tilde{a}_{r-1,r-1}^{(0)} - \tilde{a}_{r,r}^{(0)}) \tilde{a}_{r+1,r}^{(0)} \tilde{a}_{r,r+1}^{(1)} \label{eq:c1}\\
& -\tilde{a}_{r,r-1}^{(0)} \tilde{a}_{r-1,r}^{(1)} (\tilde{a}_{r+1,r+1}^{(0)} - \tilde{a}_{r,r}^{(0)})\} \{ (\tilde{a}_{r'-1,r'-1}^{(0)}- \tilde{a}_{r',r'}^{(0)}) (\tilde{a}_{r',r'}^{(1)} -c_1) (\tilde{a}_{r'+1,r'+1}^{(0)} - \tilde{a}_{r',r'}^{(0)}) \nonumber \\
& -(\tilde{a}_{r'-1,r'-1}^{(0)} - \tilde{a}_{r',r'}^{(0)}) \tilde{a}_{r'+1,r'}^{(0)} \tilde{a}_{r',r'+1}^{(1)}-\tilde{a}_{r',r'-1}^{(0)} \tilde{a}_{r'-1,r'}^{(1)} (\tilde{a}_{r'+1,r'+1}^{(0)} - \tilde{a}_{r',r'}^{(0)})\} =0, \nonumber 
\end{align}
which follows from expanding the characteristic polynomial of the matrix $\left( a_{i,j}^{{\bf a}} \right) _{i,j=0,\dots,d}$ in variable $p$ and observing the coefficient of $p^2$. 
By a direct calculation, the condition that Eq.(\ref{eq:c1}) for the variable $c_1$ has multiple roots is written as $(\alpha _3-\alpha _1)(2r+\alpha _2+\alpha _3)=0$.  If $2r+\alpha _2+\alpha _3=0$, then $r=r'$ and it contradicts that $r<r'$. By the assumption $\alpha _3 \neq \alpha _1$, it follows that Eq.(\ref{eq:c1}) for the variable $c_1$ does not have multiple roots, and the solution $E=\tilde{a}_{r,r}^{(0)} + c_1 p + \dots $ separates.

If $r+1 =r'$, then $r=- (\alpha _2+ \alpha _3+1)/2$, the eigenvalue around $\tilde{a}_{r,r}^{(0)}$ is expanded as
$E=\tilde{a}_{r,r}^{(0)} + c_{1/2} \sqrt{p} + \dots $, and $c_{1/2}$ is determined by 
\begin{align}
& c_{1/2}^2 = \tilde{a}_{r+1,r}^{(0)} \tilde{a}_{r,r+1}^{(1)}= -16(\alpha _0-\alpha _1)(\alpha _2-\alpha _3)(\alpha _0+\alpha _1-1)(\alpha _2+\alpha _3-1).
\label{eq:c1/2}
\end{align}
Since $0 \leq r+r'= -(\alpha _2+ \alpha _3) \leq 2d= -(\alpha _0+\alpha _1+\alpha _2+\alpha _3)$, we have $(1-(\alpha _0+\alpha _1))(1-(\alpha _2+\alpha _3)) > 0$.
Combining with $\alpha _2\neq \alpha _3$, it follows that, if $\alpha _0 \neq \alpha _1$, then Eq.(\ref{eq:c1/2}) for the variable $c_{1/2}$ does not have multiple roots, and the solution $E=\tilde{a}_{r,r}^{(0)} + c_{1/2} \sqrt{p} + \dots $ separates.
If $\alpha _0= \alpha _1$, then we have $a_{r+1,r}^{{\bf a}}= \tilde{a}_{r+1,r}^{(0)}=0$ for $r=- (\alpha _2+ \alpha _3+1)/2$. The eigenvalue around $\tilde{a}_{r,r}^{(0)}$ is expanded as $E=\tilde{a}_{r,r}^{(0)} +c_1 p + \dots $, and $c_{1}$ is determined by 
\begin{align}
& \{ (\tilde{a}_{r-1,r-1}^{(0)} -\tilde{a}_{r,r}^{(0)})(\tilde{a}_{r,r}^{(1)} -c_1)  -\tilde{a}_{r,r-1}^{(0)} \tilde{a}_{r-1,r}^{(1)} \} \cdot \label{eq:c1a0a1} \\
& \quad \quad \{ (\tilde{a}_{r+2,r+2}^{(0)} -\tilde{a}_{r+1,r+1}^{(0)})(\tilde{a}_{r+1,r+1}^{(1)} -c_1)  -\tilde{a}_{r+2,r+1}^{(0)} \tilde{a}_{r+1,r+2}^{(1)} \} =0 \nonumber .
\end{align}
The condition that Eq.(\ref{eq:c1a0a1}) for the variable $c_1$ has multiple roots is written as $(\alpha _2-\alpha _3)(\alpha _2+\alpha _3-1)(2\alpha _0-1)=0$. But it is impossible, because $(1-(\alpha _2+\alpha _3))(1-2\alpha _0)>0$ and $\alpha _3 \neq \alpha _2$. Hence, if $\alpha _0 = \alpha _1$, then Eq.(\ref{eq:c1a0a1}) for the variable $c_{1}$ does not have multiple roots, and the solution $E=\tilde{a}_{r,r}^{(0)} + c_{1} p + \dots $ separates.

Thus we obtain that the multiple roots $E=\tilde{a}_{r,r}^{(0)} $ at $p=0$ separates by expanding the eigenvalue as a series in $p$ or $\sqrt{p}$.

Therefore the zeros of the characteristic polynomial equation are distinct for generic periods $(2\omega _1, 2\omega _3)$.
\end{proof}

\begin{cor} (Proposition \ref{thm:dist}) \label{athm:dist}
Let $l_0,l_1,l_2,l_3$ be non-negative integers and $V$ be the vector space written as Eq.(\ref{eq:Ve}) ($l_0+l_1+l_2+l_3$: even) or Eq.(\ref{eq:Vo}) ($l_0+l_1+l_2+l_3$: odd).
We denote the monic characteristic polynomial of the operator $H^{(l_0,l_1,l_2, l_3)}$ on the space $V$ by $P(E)$.
Then the roots of the equation $P(E)=0$ are distinct for generic periods $(2\omega _1, 2\omega _3)$.
\end{cor}
\begin{proof}
It is shown in the proof of \cite[Theorem 3.2]{Tak1} (or \cite[Proposition 3.9]{TakP}) that any two characteristic polynomials of distinct subspaces in Eq.(\ref{eq:Ve0}) or Eq.(\ref{eq:Vo0}) does not have common roots. Hence it is sufficient to show that the characteristic polynomial of the operator $H^{(l_0,l_1,l_2, l_3)}$ on any space listed in Eq.(\ref{eq:Ve0}) or Eq.(\ref{eq:Vo0}) does not have multiple zeros for generic periods $(2\omega _1, 2\omega_3)$. In Proposition \ref{aprop:distrt}, it is shown that, if $\alpha _i \neq \alpha_j$ for some $i,j \in \{0,1,2,3\}$, then the characteristic polynomial does not have multiple zeros for generic periods.
In Eqs.(\ref{eq:Ve0}, \ref{eq:Vo0}), the case $\alpha _0 =\alpha _1 =\alpha _2 =\alpha _3 $ for $V _{\alpha _0, \alpha _1, \alpha _2, \alpha _3}$ appears only for the case $l_0=l_1=l_2=l_3$ and $\alpha _0 =\alpha _1 =\alpha _2 =\alpha _3 =-l_0$. We set $g=l_0 (=l_1=l_2=l_3)$.

For the case $l_0=l_1=l_2=l_3=g$, the operator $H^{(g,g,g,g)}$ is expressed as 
\begin{equation*}
H^{(g,g,g,g)}= -\frac{d^2}{dx^2} + 4g(g+1)\wp (2x),
\end{equation*}
and the finite-dimensional space $V(=V_{-g, -g, -g, -g})$  for the case $l_0=l_1=l_2=l_3=g$ coincides with the space $V$ for the case $l_0=g$ and $l_1=l_2=l_3=0$ by replacing basic periods $(2\omega _1, 2\omega _3) \rightarrow (\omega _1, \omega _3)$. For the case $l_0 \neq 0$ and $l_1=l_2=l_3=0$, the corresponding proposition is proved in Proposition \ref{aprop:distrt} or \cite[\S 23.4]{WW}. 
Thus Corollary \ref{athm:dist} is proved.
\end{proof}

\end{document}